\documentclass[10pt,a4paper]{amsart} 
\usepackage{multicol}

\usepackage[utf8]{inputenc}
\usepackage[T1]{fontenc}

\usepackage[colorlinks=true, citecolor=purple, urlcolor=blue, linkcolor=purple%, pagebackref
]{hyperref}
\usepackage[noabbrev]{cleveref}
\usepackage{amsfonts}
\usepackage{amssymb}
\usepackage{amsthm} 

\usepackage[all]{xy}

\usepackage{mathtools} % coloneqq

\usepackage{mathdots} % iddots

\usepackage{mathrsfs}

\usepackage{dutchcal} % new calligraphic font with small letters

\usepackage{tikz}
\usetikzlibrary{cd}

\numberwithin{equation}{subsection}

\newtheorem{main}{Theorem}
\crefname{main}{theorem}{theorems}
\Crefname{main}{Theorem}{Theorems}

\newtheorem*{conjecture*}{Conjecture}

\newtheorem{theorem}[equation]{Theorem}

\newtheorem{proposition}[equation]{Proposition}
\newtheorem*{proposition*}{Proposition}

\newtheorem*{corollary*}{Corollary}
\newtheorem{lemma}[equation]{Lemma}

\theoremstyle{definition}

\theoremstyle{remark}
\newtheorem{remark}[equation]{Remark}
\newtheorem*{remark*}{Remark}
\newtheorem{example}[equation]{Example}
\newtheorem*{example*}{Example}

\newtheorem*{examples*}{Examples}

\newcommand{\hooklongrightarrow}{\lhook\joinrel\longrightarrow}

\DeclareMathOperator{\id}{id}
\newcommand{\C}{\mathbf{C}}

\newcommand{\Z}{\mathbf{Z}}

\DeclareMathOperator{\gr}{gr}
\DeclareMathOperator{\Spec}{Spec}
\DeclareMathOperator{\Sym}{Sym}
\newcommand{\Hgy}{\mathrm{H}}
\newcommand{\Cpx}{\mathrm{C}}

\newcommand{\git}{\big/\!\!\big/}

\DeclareMathOperator{\ad}{ad}

\newcommand{\kfk}{\mathfrak k}
\newcommand{\Ga}{\mathbf{G}_\mathrm{a}}
\newcommand{\Gm}{\mathbf{G}_\mathrm{m}}

\newcommand{\g}{\mathfrak g}
\newcommand{\n}{\mathfrak n}

\renewcommand{\sl}{\mathfrak{sl}}

\newcommand{\F}{\mathrm F}
\newcommand{\J}{\mathrm{J}}
\newcommand{\U}{\mathrm U}
\newcommand{\V}{\mathrm V}
\newcommand{\W}{\mathrm W}

\newcommand\NO[1]{\mathopen{\boldsymbol{\vcentcolon}}#1\mathclose{\boldsymbol{\vcentcolon}}}

\title{A decomposition theorem for minimal W-algebras in type A}

\author{Thibault Juillard}
\address{TJ: Universität Hamburg, Hamburg, Germany}
\email{thibault.juillard@uni-hamburg.de}

\usepackage[backend=biber,style=alphabetic,autocite=plain, maxbibnames=10]{biblatex}
\begin{document}
	
\begin{abstract}
	To any simple Lie algebra, one can associate a minimal affine W-algebra. It is a vertex algebra which is obtained from the universal affine vertex algebra by a chiral Hamiltonian reduction procedure. This reduction is determined by data coming from the minimal nilpotent orbit. It is a chiral quantisation of the Slodowy slice associated to this minimal orbit.
	
	In this paper, when the simple Lie algebra is of type A, we prove that the minimal W-algebra is isomorphic to a commutant subalgebra in a localisation of the universal affine vertex algebra. Then, we establish a decomposition theorem \emph{à la} Losev: the localised universal affine vertex algebra decomposes as the tensor product of the minimal affine W-algebra and another vertex subalgebra. As a direct corollary, we retrieve the celebrated result of inverse Hamiltonian reduction for the minimal W-algebra.
	
	We also prove an analogue of this result for the corresponding finite W-algebra. In both cases, our approach relies on the geometry of the associated Slodowy slice.
\end{abstract}

\maketitle

\tableofcontents

\section{Introduction}

\subsection{Decomposition theorem and inverse reduction}

To a simple Lie algebra $\g$ and a nilpotent element $f$, one can associate two interesting algebraic structures: 
\begin{enumerate}
	\item The \emph{finite W-algebra} $\U(\g, f)$: it is an associative algebra, constructed as the quantum Hamiltonian reduction of the universal enveloping algebra $\U(\g)$ \cite{kostant1978whittaker,premet2002special,gan2002quantization}.
	
	\item The \emph{affine W-algebra} $\W^k(\g, f)$: it is a vertex algebra, constructed as the chiral Hamiltonian reduction of the universal affine vertex algebra $\V^k(\g)$ \cite{feigin1990quantization,kac2003quantum,kac2004quantum}.
\end{enumerate}
The finite W-algebra is related to the affine W-algebra: the former is the associated \emph{Zhu algebra} of the latter \cite{desole2006finite}. Therefore, both their representation theories are related \cite{zhu1996modular}.

Losev introduced and studied an $\hbar$-adic analogue $\U(\g, f)_\hbar$ of the finite W-algebra \cite{losev2010quantized}, which is a complete $\C[\![\hbar]\!]$-algebra. He also established his celebrated \emph{decomposition theorem}, that is a $\hbar$-adic algebra isomorphism
\begin{equation} \label{equation:decomposition-theorem}
	\U(\g)_\hbar \cong A_\hbar \widehat\otimes_{\C[\![\hbar]\!]} \U(\g, f)_\hbar, 
\end{equation}
where $\U(\g)_\hbar$ is an $\hbar$-adic analogue of the enveloping algebra and $A_\hbar$ is a $\hbar$-adic Weyl algebra. Such a result is a crucial tool to connect the representation theories of $\g$ and the finite W-algebra $\U(\g, f)$. His proof is based on the geometry of the adjoint orbit of $f$ and its transverse slice $S_f$, called the \emph{Slodowy slice} of $f$.

For affine W-algebras, a similar kind of phenomenon has been intensively studied under the name of the \emph{inverse Hamiltonian reduction}. By definition, inverse reduction holds if there exists a vertex algebra embedding of the form
\begin{equation} \label{equation:inverse-reduction}
	\V^k(\g) \hooklongrightarrow  \Pi_f \otimes \W^k(\g, f), 
\end{equation}
where $\Pi_f$ is a certain localisation of a vertex algebraic analogue of a Weyl algebra. Such an embedding has been used to study weight modules of the simple quotient of $\V^k(\g)$ by using the representation theory of $\Pi_f$ and $\W^k(\g, f)$ \cite{adamovic2019realizations,adamovic2021realisation,adamovic2024relaxed,fasquel2025modularity}. 

Such inverse reductions have been established following different approaches:
\begin{itemize}
	\item When $\g$ is of small rank, inverse reduction can be constructed by using explicit generator and relations descriptions of W-algebras \cite{semikhatov1994inverting,adamovic2019realizations,adamovic2021realisation,adamovic2024relaxed}.
	
	\item By using the screening operator description of W-algebras \cite{genra2020screening}, inverse reduction has been established for any nilpotent $f$ when $\g$ is of type A \cite{fehily2023subregular,fehily2024inverse,fasquel2026connecting,fasquel2026quantum}.
\end{itemize}

To prove the decomposition theorem \eqref{equation:decomposition-theorem}, Losev uses the fact that the finite W-algebra $\U(\g, f)$ (and its $\hbar$-adic analogue) is a quantisation of the Poisson algebra $\C[S_f]$ of polynomial functions on the Slodowy slice $S_f$ \cite{premet2002special,gan2002quantization}. Such a relation also holds for the affine W-algebra: $\W^k(\g, f)$ carries a filtration \cite{li2004vertex,li2005abelianizing} such that the associated graded vertex algebra $\gr \W^k(\g, f)$ is commutative and isomorphic to $\C[\J_\infty S_f]$, the algebra of polynomial functions on the \emph{arc space} $\J_\infty S_f$ of the Slodowy slice \cite{desole2006finite,arakawa2015associated}. However, the approaches described above to establish inverse reduction \eqref{equation:inverse-reduction} do not make use of the Slodowy slice.

The Poisson structure of the Slodowy slice can be described by Hamiltonian reduction \cite{kostant1978whittaker,premet2002special,gan2002quantization}. There is a unipotent group $G_{\geqslant 1}$ that acts on the dual space $\g^\vee$, and a moment map $\mu \colon \g^\vee \rightarrow (\g_{\geqslant 1})^\vee$, such that $S_f$ is isomorphic to the Hamiltonian reduction $\mu^{-1}(0) / G_{\geqslant 1}$. In fact, there is a $G_{\geqslant 1}$-equivariant isomorphism
\begin{equation} \label{equation:gan-ginzburg} 
	G_{\geqslant 1} \times S_f \cong \mu^{-1}(0) 
\end{equation}
of varieties (see \Cref{theorem:gan-ginzburg-isomorphism}). 

This isomorphism is similar to the decomposition theorem \eqref{equation:decomposition-theorem} and the inverse reduction \eqref{equation:inverse-reduction}, with $S_f$ playing the role of the corresponding W-algebras. However, it is not the right geometric analogue because the dimension of the two varieties is too small and these varieties are not Poisson. The right analogue of \eqref{equation:decomposition-theorem} and \eqref{equation:inverse-reduction} must involve Poisson varieties of the same dimension as the Lie algebra $\g$. 

\subsection{A geometric approach: the affine Grassmannian}

Recently, a geometrical understanding of inverse reduction has been developed by using the affine Grassmannian \cite{kamnitzer2022hamiltonian,butson2025inverse,butson2025inverse-affine}. By the Mirkovi\'c--Vybornov correspondence \cite{mirkovic2022comparison}, Slodowy slices in type A are in correspondence with slices in the affine Grassmannian. A quantum analogue of this statement relates finite W-algebras in type A and truncated shifted Yangians \cite{webster2020quantum}. Such a correspondence is not established yet for affine W-algebras.

In \cite{kamnitzer2022hamiltonian,butson2025inverse}, the multiplication map of the affine Grassmannian is used to provide an inverse reduction statement for Slodowy slices and finite W-algebras, through the Mirkovi\'c--Vybornov correspondence. In particular, the multiplication map gives an analogue of the isomorphism \eqref{equation:gan-ginzburg} with the correct dimension: one side is an open subset of $\mathrm{T}^\vee G_{\geqslant 1} \times S_f$ (the product of the cotangent bundle of the unipotent group with the Slodowy slice), and the other side is an open subset of $\g^\vee$.

In \cite{butson2025inverse-affine}, Butson and Nair use the chiral quantisation theory they have developed in \cite{butson2026deformation} and the inverse reduction for Slodowy slices to establish inverse reduction \eqref{equation:inverse-reduction} in type A \cite[Section 6.3]{butson2025inverse-affine}. Note that only the existence of such an isomorphism is provided by the machinery of \cite{butson2026deformation}, but there is no explicit description of the map. The rough idea is to prove that both sides of the isomorphism are chiral quantisation of the same Poisson variety, and then use the classification of such quantisations to conclude that they must be isomorphic as vertex algebras. This argument is the vertex algebraic counterpart to Losev's approach for finite W-algebras \cite{losev2010quantized}.

\subsection{Our contribution}

In the present paper, we propose a new approach to inverse reduction that starts from the construction of the affine W-algebra as a chiral Hamiltonian reduction \cite{feigin1990quantization,kac2003quantum}. This reduction is concretely achieved by using the \emph{BRST (Becchi–Rouet–Stora–Tyutin) cohomology} functor (\Cref{section:brst}) $\Hgy^0_f$,
\[ \W^k(\g, f) \coloneqq \Hgy^0_f(\V^k(\g)). \]
Our construction is achieved in type A, for the minimal nilpotent orbit. For the rest of the paper, let $\g \coloneqq \sl_d$ be a simple Lie algebra of type A, $d \geqslant 2$ being an integer. Let $(e, h, f)$ be a minimal $\sl_2$-triple in $\g$.

The starting point of our construction is to establish an analogue of \eqref{equation:gan-ginzburg} between varieties of the correct dimension, without using the geometry of the affine Grassmannian. The $\sl_2$-triple equips $\g$ with a grading and $G_{\geqslant 1}$ is the unipotent group associated to the positive degree subalgebra $\g_{\geqslant 1}$. In \Cref{subsection:slodowy-quotient-open}, we construct a bigger group $K \ltimes G_{\geqslant 1}$, and its associated Lie algebra $\kfk \ltimes \g_{\geqslant 1}$ (see \Cref{lemma:afk}), such that the action of this group provides the desired isomorphism.

\begin{main}[\Cref{proposition:dilation-gan-ginzburg} and \Cref{theorem:decomposition-slodowy}] \label{main:geometric-setting}
	The group $K \ltimes G_{\geqslant 1}$ acts on the distinguished open set 
	\[ \g^\vee[e^{-1}] \coloneqq \Spec \big((\Sym \g)[e^{-1}]\big), \] 
	such that there is a $K \ltimes G_{\geqslant 1}$-equivariant isomorphism
	\[ (K \ltimes G_{\geqslant 1}) \times S_f \cong \g^\vee[e^{-1}] \]
	induced by the inclusion of the minimal Slodowy slice $S_f$ in $\g^\vee[e^{-1}]$.
\end{main}

By use of a ``quantisation commutes with reduction'' argument (\Cref{proposition:reduction-quantisation}), we can use \Cref{main:geometric-setting} to prove the main result of this paper.

\begin{main}[\Cref{theorem:affine-w-algebra-coset,theorem:decomposition-chiral}] \label{main:chiral-setting}
	One can construct a \emph{localisation} $\V^k(\g)[e^{-1}]$ of the universal affine vertex algebra $\V^k(\g)$ such that:
	\begin{enumerate}
		\item The affine W-algebra $\W^k(\g, f)$ is naturally isomorphic to a commutant subalgebra of $\V^k(\g)[e^{-1}]$:
		\[ \W^k(\g, f) \cong \{a \in \V^k(\g)[e^{-1}] \quad | \quad \text{$a$ commutes with $\V^k(\kfk \ltimes \g_{\geqslant 1})$} \}, \]
		where $\V^k(\kfk \ltimes \g_{\geqslant 1})$ is the vertex subalgebra of $\V^k(\g)$ generated by $\kfk \ltimes \g_{\geqslant 1}$.
		
		\item There is a vertex algebra isomorphism
		\begin{equation} \label{equation:main-decomposition-theorem}
			\V^k(\kfk \ltimes \g_{\geqslant 1})[e^{-1}] \otimes \W^k(\g, f) \overset{\sim}{\longrightarrow}  \V^k(\g)[e^{-1}] 
		\end{equation}
		given by the inclusion of the two mutually commuting vertex subalgebras $\V^k(\kfk \ltimes \g_{\geqslant 1})[e^{-1}]$ and $\W^k(\g, f)$ in $\V^k(\g)[e^{-1}]$.
	\end{enumerate}
\end{main}

The isomorphism \eqref{equation:main-decomposition-theorem} is an analogue of Losev's decomposition theorem, for non-$\hbar$-adic affine W-algebras. The theorem holds for any level $k$, no genericness is assumed. Note that, unlike the construction of \cite{butson2025inverse-affine}, the isomorphism that we construct is given by an explicit map. 

By \Cref{theorem:commutant-localisation-weyl}, the vertex subalgebra $	\V^k(\kfk \ltimes \g_{\geqslant 1})[e^{-1}]$ is in fact the localised Weyl vertex algebra $\Pi_f$ that appears in the literature on inverse reduction. \emph{Therefore, inverse reduction follows from our theorem in a straightforward way.} Indeed, localisation gives a vertex algebra embedding $\V^k(\g) \hookrightarrow \V^k(\g)[e^{-1}]$ (\Cref{subsection:localisation-vertex}). By composing this embedding with the inverse of the decomposition isomorphism \eqref{equation:main-decomposition-theorem}, we get the inverse reduction \eqref{equation:inverse-reduction}. 

The vertex algebra homomorphisms appearing in \Cref{main:chiral-setting} are not difficult to construct. The crucial point is to prove that they are isomorphisms. To do so, we use the increasing Li filtrations on vertex algebras (\Cref{subsection:vertex-algebras}). Taking the associated graded objects provide commutative algebras that correspond to polynomial functions on schemes:
\[ \gr \V^k(\g) \cong \C[\J_\infty \g^\vee],  \quad \gr \big(\V^k(\g)[e^{-1}]\big) \cong \C[\J_\infty \g^\vee][e^{-1}], \quad   \gr \W^k(\g, f) \cong \C[\J_\infty S_f], \]
where $\J_\infty X$ denotes the arc scheme of a scheme $X$ (\Cref{subsection:arc-schemes}).

Finally, \Cref{main:geometric-setting} also leads to a decomposition theorem and inverse reduction statement for the minimal finite W-algebra.

\begin{main}[\Cref{theorem:finite-w-algebra-coset,theorem:decomposition-quantum}] \label{main:quantum-setting}
	The minimal finite W-algebra $\U(\g, f)$ is naturally isomorphic to the commutant subalgebra
	\[ (\U(\g)[e^{-1}])^{\kfk \ltimes \g_{\geqslant 1}} \coloneqq \{a \in \U(\g)[e^{-1}] \quad | \quad \text{$a$ commutes with $\kfk \ltimes \g_{\geqslant 1}$} \}. \]
	Moreover, there is an algebra isomorphism
	\begin{equation} \label{equation:main-decomposition-theorem-finite}
		\U(\kfk \ltimes \g_{\geqslant 1})[e^{-1}] \otimes \U(\g, f) \overset{\sim}{\longrightarrow}  \U(\g)[e^{-1}]. 
	\end{equation}
\end{main}

\subsection{The $\sl_2$ example}

We illustrate our result in the easiest case, when $\g = \sl_2$ and
\[ e = \left(\begin{smallmatrix}
	0 & 1 \\
	0 & 0
\end{smallmatrix}\right), \quad h = \left(\begin{smallmatrix}
	1 & 0 \\
	0 & -1
\end{smallmatrix}\right) \quad \text{and} \quad f = \left(\begin{smallmatrix}
	0 & 0 \\
	1 & 0
\end{smallmatrix}\right). \]

\subsubsection{The geometric setting}

The symmetric algebra $\Sym \sl_2$ has a Poisson structure given by
\[ \{h, e\} = 2 e, \quad \{e, f\} = h \quad \text{and} \quad \{h, f\} = - 2 f. \]

In this case the Slodowy slice is given by the explicit description:
\[ S_f \coloneqq \Big\{\begin{pmatrix}
	0 & \star \\
	1 & 0
\end{pmatrix}\Big\}. \]
By the work of Kostant \cite{kostant1978whittaker}, the Poisson centre $\mathrm{Z}(\Sym \sl_2)$ is generated by the Casimir element
\[ c \coloneqq \frac{1}{2} h^2 + 2 e f \] 
and isomorphic to $\C[S_f]$.

There is an obvious Poisson algebra homomorphism
\[ \Sym(\C h \oplus \C e) \otimes \mathrm{Z}(\Sym \sl_2) \longrightarrow \Sym \sl_2. \]
Note that if we allow $e$ to be invertible, then we can write
\[ f = e^{-1} \frac{1}{2}(c - \frac{1}{2} h^2). \]
Therefore we have a Poisson algebra isomorphism
\[ \Sym(\C h \oplus \C e)[e^{-1}] \otimes \mathrm{Z}(\Sym \sl_2) \longrightarrow (\Sym \sl_2)[e^{-1}]. \]

\subsubsection{The quantum setting}

The enveloping algebra $\U(\sl_2)$ has commutations relations given by
\[ [h, e] = 2 e, \quad [e, f] = h \quad \text{and} \quad [h, f] = - 2 f. \]

By the work of Kostant \cite{kostant1978whittaker}, the finite W-algebra $\U(\sl_2, f)$ is isomorphic to the centre $\mathrm{Z}(\U(\sl_2))$, that is generated by the Casimir element
\[ \widetilde c \coloneqq \frac{1}{2} h^2 + 2 e f - h. \] 

There is an obvious algebra homomorphism
\[ \U(\C h \oplus \C e) \otimes \mathrm{Z}(\U(\sl_2)) \longrightarrow \U(\sl_2). \]
Note that if we allow $e$ to be invertible, then we can write
\[ f = e^{-1} \frac{1}{2}(\widetilde{c} - \frac{1}{2} h^2 + h). \]
Therefore we have an algebra isomorphism
\[ \U(\C h \oplus \C e)[e^{-1}] \otimes \mathrm{Z}(\U(\sl_2)) \cong \U(\sl_2)[e^{-1}]. \]

\subsubsection{The chiral setting at critical level}

The universal affine vertex algebra $\V^k(\sl_2)$ has operator product expansions given by
\[ \begin{multlined}
	h(z_1) h(z_2) \sim \frac{1}{(z_1 - z_2)^2} 2 k \id, \quad h(z_1) e(z_2) \sim \frac{1}{z_1 - z_2} 2 e(z_2) , \\ e(z_1) f(z_2) \sim \frac{1}{z_1 - z_2} h(z_2) + \frac{1}{(z_1 - z_2)^2} k \id \\ \text{and} \quad h(z_1) f(z_2) \sim - \frac{1}{z_1 - z_2} 2 f(z_2).
\end{multlined} \]
To simplify the exposition, we consider only the critical level $k = -2$. For the general case, see \Cref{section:sl2-example}.

By the work of Feigin--Frenkel \cite{feigin1992affine}, at the critical level, the affine W-algebra $\W^{-2}(\sl_2, f)$ is isomorphic to the centre of $\V^{-2}(\sl_2)$ and is generated by
\[ \widehat c \coloneqq \frac{1}{2} \NO{h h} + \NO{e f} + \NO{f e} = \frac{1}{2} \NO{h h} + 2 \NO{e f} - \partial h, \]
where $\partial$ denotes the translation operator and $\NO{}$ denotes the normally ordered product. 

There is an obvious vertex algebra homomorphism
\[ \V^{-2}(\C h \oplus \C e) \otimes \W^{-2}(\sl_2, f) \longrightarrow \V^{-2}(\sl_2). \]
Note that if we allow $e$ to be invertible, then we can write
\[ \NO{e^{-1} (\NO{e f})} = \NO{e^{-1} \frac{1}{2}(\widehat{c} - \frac{1}{2} \NO{h h} + \partial h)}. \]
The non-associativity of the normally ordered product is controlled by the following Borcherds identity:
\[ \begin{multlined}
	f = {(\NO{e^{-1} e}) f} = \NO{e^{-1} (\NO{e f})} + P, 
	\\ \text{where} \quad P \coloneqq \sum_{i \geqslant 0} \big(e^{-1}{}_{(-2-i)}(e_{(i)} f) + e{}_{(-2-i)}(e^{-1}{}_{(i)} f)\big), 
\end{multlined}  \]
with the notation $\bullet_{(i)}\bullet$ for the $i$-products in the vertex algebra.

The operator product expansion of $e$ and $f$ says that $e_{(i)}f$ belongs to the subalgebra $\V^{-2}(\C h \oplus \C e)$ for any nonnegative integer $i$. Then, \cite[Formula (3.1)]{arakawa2026faithful} implies that $e^{-1}{}_{(i)}f$ belongs to $\V^{-2}(\C h \oplus \C e)[e^{-1}]$ for any nonnegative integer $i$. Therefore, $P$ belongs to $\V^{-2}(\C h \oplus \C e)[e^{-1}]$ and
\[ f = \NO{e^{-1} \frac{1}{2}(\widehat{c} - \frac{1}{2} \NO{h h} + \partial h)} + P. \]
See \Cref{proposition:computation-f} for an explicit formula.

Therefore we get a vertex algebra isomorphism
\[ \V^{-2}(\C h \oplus \C e)[e^{-1}] \otimes \W^{-2}(\sl_2, f) \cong \V^{-2}(\sl_2)[e^{-1}]. \]

\subsection{Future work and applications}

\subsubsection{Reduction by stages in type A}

In a work in progress \cite{juillard2027}, we extend our construction to inverse reduction by stages in type A. By using our previous work with Genra \cite{genra2024reduction,genra2025reduction}, for two adjacent hook-type nilpotent elements $f_1$ and $f_2$ in $\g = \sl_d$ such that the orbit of $f_2$ is bigger than the orbit of $f_1$, we construct an explicit vertex algebra isomorphism
\[ \Pi_{f_1, f_2} \otimes \W^k(\g, f_2) \cong \W^k(\g, f_1)[e_0{}^{-1}], \]
where $e_0$ is a localisable element in $\W^k(\g, f_1)$ and $\Pi_{f_1, f_2}$ is a localised Weyl vertex algebra. By induction, we can deduce a decomposition theorem for any hook-type nilpotent element $f$ in $\g$:
\[ \Pi_{f} \otimes \W^k(\g, f) \cong \V^k(\g)[S], \]
where $S$ is a localisable subset in $\V^k(\g)$ and $\Pi_{f}$ is a localised Weyl vertex algebra.

\subsubsection{Minimal W-algebra in any type}

If $\g$ is not of type A, \Cref{main:geometric-setting} no longer holds: we need to replace the open set $\g^\vee[e^{-1}]$ by a suitable étale cover. On the noncommutative side, this means that we need to replace the localisation embedding $\V^k(\g) \hookrightarrow \V^k(\g)[e^{-1}]$ by a more complicated extension of $\V^k(\g)$ to get the right analogue of \Cref{main:chiral-setting} for a general $\g$. 

With Christopher Raymond, we aim to establish a decomposition theorem for any simple Lie algebra $\g$ of rank $2$, that is to say of type $\mathrm{A}_2$, $\mathrm{B}_2$ or $\mathrm{G}_2$. Note that an inverse reduction in type $\mathrm{A}_2$ is already known and is the one retrieved by \Cref{main:chiral-setting}. However, replacing $\V^k(\sl_3) \hookrightarrow \V^k(\sl_3)[e^{-1}]$ by a more complicated extension of $\V^k(\sl_3)$ should lead to another inverse reduction. 

As mentioned earlier, inverse reduction has proved useful to study weight modules over the simple quotient of $\V^k(\g)$ \cite{adamovic2019realizations,adamovic2021realisation,adamovic2024relaxed,fasquel2025modularity}. We want to pursue these studies in these new cases. 

\subsection{Organisation of the paper}

In \Cref{section:geometric-setting}, we recall the construction of the Slodowy slice as Hamiltonian reduction and establish the decomposition theorem for the minimal Slodowy slice (\Cref{theorem:decomposition-slodowy}). An important result for the rest of the paper is \Cref{proposition:dilation-gan-ginzburg}, that is an analogue of the Gan--Ginzburg isomorphism (\Cref{theorem:gan-ginzburg-isomorphism}) for an open set in $\g^\vee$. 

In \Cref{section:localisation-filtration}, we collect technical results on filtrations on noncommutative algebras and localisation in the noncommutative setting. We also give some results to compute invariant subspaces of filtered modules for an algebraic group.

The results of the two previous sections are applied in \Cref{section:quantum-setting} to prove that the minimal finite W-algebra is an invariant subalgebra of the localised enveloping algebra (\Cref{theorem:finite-w-algebra-coset}). Then, we establish the decomposition theorem for the minimal finite W-algebra (\Cref{theorem:decomposition-quantum}).

In \Cref{section:vertex-algebra}, we collect results about infinite type objects: vertex algebras and arc spaces. We also recall the Borisov localisation theory of vertex algebras and its compatibility with the Li filtration. We give some technical results on computing invariant subspaces of filtered modules for the projective limit of a system of algebraic groups. 

In \Cref{section:brst}, we recall the BRST construction of chiral Hamiltonian reduction and we give a sufficient condition to prove that ``quantisation commutes with reduction'' (\Cref{proposition:reduction-quantisation}).

In \Cref{section:chiral-setting}, we prove that the minimal affine W-algebra is an invariant subalgebra of the localised universal vertex algebra (\Cref{theorem:affine-w-algebra-coset}). Then, we establish the decomposition theorem for the minimal affine W-algebra (\Cref{theorem:decomposition-chiral}). Finally, we prove that the subalgebra of the localisation commuting with the W-algebra is a localised Weyl vertex algebra (\Cref{theorem:commutant-localisation-weyl}).

In \Cref{section:sl2-example}, we give the example of our construction, in the vertex-algebraic setting, for $\g = \sl_2$. In particular, we make the computations at any level.

\subsection{Notations and conventions}

Unless otherwise stated, all objects (vector spaces, algebras, schemes, tensor products, ...) are defined over the field of complex numbers $\C$. Denote by $\Z$ the ring of integers. Denote by $\delta_{a=b}$ the Kronecker symbol that is equal to $1$ if $a=b$ holds, $0$ otherwise. 

If $V$ is a vector space, denote by $\Sym V$ its symmetric algebra, $\operatorname{Alt} V$ its exterior algebra and $V^\vee$ its dual vector space.

For any affine scheme $X = \Spec R$, denote by $\C[X] \coloneqq R$ its coordinate ring. If $\phi : X \rightarrow Y$ is a homomorphism of affine schemes, the comorphism is denoted by~$\phi^* : \C[Y] \rightarrow \C[X]$. If $R$ is a $\C$-algebra and $X$ is a scheme, then $X(R)$ denotes the set of $R$-points of $X$, that is to say the set of unital $\C$-algebra homomorphisms from $\C[X]$ to $R$.

If $S$ is a multiplicative subset in $R$, denote by $R[S^{-1}]$ the corresponding localisation. Denote by $X[S^{-1}] \coloneqq \Spec R[S^{-1}]$ the corresponding open affine subscheme. If $S$ is generated by one element $F$ in $R$, then we adopt the notation $R[F^{-1}] \coloneqq R[S^{-1}]$ and $X[F^{-1}] \coloneqq \Spec R[F^{-1}]$.  

If $X$ is a variety (reduced affine scheme of finite type), the notation $x \in X$ means that $x$ is a $\C$-point of $X$.  A linear algebraic group is an affine variety with a group scheme structure. The multiplicative group is denoted by~$\Gm \coloneqq \Spec \C[t, t^{-1}]$ and is identified with its $\C$-point group of invertible complex numbers $\C^\times$. The additive group is denoted by~$\Ga \coloneqq \Spec \C[t]$ and is identified with its $\C$-point group of complex numbers $\C$.

If $G$ is a linear algebraic group, its Lie algebra is denoted by $\g$. Actions of an algebraic group on a scheme or a vector space are assumed to be algebraic, that is to say given by a co-action \cite{milne2017algebraic}. 

If both $G$ and $\g$ act on a vector space $V$, for any vector $v$ in $V$ set
$$ G^v \coloneqq \{g \in G \quad | \quad g \cdot v = v \} \quad \text{and} \quad \g^v \coloneqq \{x \in \g \quad | \quad x \cdot v = 0 \}. $$
Note that $\g^v$ is the Lie algebra of $G^v$. The subspaces of invariants are denoted by
\[ \begin{aligned}
	V^G & \coloneqq \{v \in V \quad | \quad \text{for all $g \in G$},~ g \cdot v = v \} \\
	\text{and} \quad V^\g & \coloneqq \{v \in V \quad | \quad \text{for all $x \in \g$},~ x \cdot v = 0 \}
\end{aligned}\]

If $G$ acts on an affine scheme $X \coloneqq \Spec R$, the affine categorical quotient is the affine scheme
\[ X \git G \coloneqq \Spec R^G. \]

\subsection{Acknowledgements}

The idea of this project arose in a discussion with Naoki Genra on the $\g = \sl_2$ case. I am very grateful to him for many years of stimulating discussions and collaboration. For the same reason, I am grateful to Christopher Raymond. In particular, our conversations were essential to achieve the proof of \Cref{theorem:commutant-localisation-weyl} and the computations of \Cref{section:sl2-example}.

I thank Anne Moreau and Sven M\"oller for their advice during the preparation of this paper and for their careful proofreading. I thank Tomoyuki Arakawa, Gwyn Bellamy, Thomas Creutzig and David Ridout for enlightening discussions. I also thank Tomoyuki Arakawa for communicating the results of \cite{arakawa2026faithful} to me. I thank Felipe Ruiz and Tanguy Vernet for their feedback on the introduction.

The $\g = \sl_2$ case has been computed independently by Tomoyuki Arakawa and Thomas Creutzig. 

This research was supported through the Emmy Noether Programme by
the \emph{Deutsche Forschungsgemeinschaft} (project number 460925688).

\subsection{Publication overlap}

During the proofreading of the present paper, we noticed the prepublication of Genqiang Liu \cite{liu2026root}. The main result of Liu is the same as our \Cref{main:quantum-setting}, but for $\g$ being of any type. 

\subsection{AI disclosure}

No generative or agentive AI model was used to carry out this research. A first version of this paper was written without AI assistance. \emph{DeepSeek~V4.1~Flash} was used by the author to produce a report on grammar, reference and notation mistakes. All the changes were then implemented manually by the author to produce the current version of the paper.

\section{Decomposition theorem in the Poisson geometric setting}

\label{section:geometric-setting}

\subsection{The coadjoint module for a Lie algebra}

Let $d \geqslant 2$ be an integer. In this paper, $\g$ will denote the Lie algebra $\sl_d$ of square matrices of size $d$ whose trace equals $0$. Denote by $x_{i,j}$ the elementary matrix whose coefficients are all $0$, except in position $(i,j)$ where it is one. We identify the dual vector space $\g^\vee$ with the affine Poisson variety
\[ \g^\vee \coloneqq \Spec (\Sym \g), \]
where the Poisson bracket on $\Sym \g$ extends the Lie bracket on $\g$.

Denote by $G$ the adjoint group of $\g$: it is the connected component of the identity in the group of Lie algebra automorphisms of $\g$. Explicitly, it is isomorphic to the projective linear group $\mathrm{PGL}_d$. Denote by $\cdot$ the adjoint action of $G$ on $\g$ and the coadjoint action on $\g^\vee$. The action on $\g^\vee$ is by Poisson automorphisms.

Denote by 
\[ (x \mid y) \coloneq \operatorname{trace}(xy), \quad \text{for} \quad x, y \in \g,  \]
the normalised Killing form on $\g$. It is $G$-invariant and non-degenerate, therefore it induces a $G$-equivariant isomorphism
\[ \kappa \colon \g \overset{\sim}{\longrightarrow} \g^\vee, \quad x \longmapsto (x \mid \bullet). \]

\subsection{The Hamiltonian action associated with the minimal nilpotent orbit}

We recall the classical construction of the Poisson structure on the minimal Slodowy slice. It is a particular case of a more general construction for any nilpotent orbit \cite{kostant1978whittaker,lynch1979generalized,premet2002special,gan2002quantization,losev2010quantized}. 

Consider the minimal $\sl_2$-triple of $\g$ given by
\[ e \coloneqq x_{1,2}, \quad h \coloneqq x_{1,1} - x_{2,2} \quad \text{and} \quad f \coloneq x_{2,1}. \]
We consider the diagonal element $H \coloneqq x_{1,1} - \frac{1}{d} \sum_{i=1}^{d} x_{i,i}$. Because it is semisimple, its adjoint action defines a Lie algebra grading on $\g$. It is explicitly given by:
\[ \g_{-1} = \bigoplus_{i=2}^d \C x_{i,1}, \quad \g_0 = \bigoplus_{i,j=2}^d \C (x_{i,j} - \delta_{i=j}x_{1,1}), \quad \g_{1}  = \bigoplus_{j=2}^d \C x_{1,j}. \]

It is an even good $\frac{1}{2} \Z$-grading in the sense of \cite{elashvili2005classification,brundan2007good}. In particular, the adjoint action of $f$ gives an injective map $\g_1 \hookrightarrow \g_0$ and a surjective map $\g_0 \twoheadrightarrow \g_{-1}$.

We denote the dual grading by
\[ \g^\vee = (\g^\vee)_1 \oplus (\g^\vee)_0 \oplus (\g^\vee)_{-1}, \]
where the subspace $(\g^\vee)_i \subseteq \g^\vee$ is defined as the image of the dual of the projection $\g \twoheadrightarrow \g_{-i}$ for all $i$ in $\Z$. In particular $(\g^\vee)_i \cong (\g_{-i})^\vee$.

The positive part of the grading $\g_{\geqslant 1} = \g_1$ is an Abelian nilpotent subalgebra of $\g$. It is the Lie algebra of the following unipotent subgroup of $G$:
\[ G_{\geqslant 1} \coloneqq \bigg\{\left(\begin{smallmatrix}
	1 & \star & \cdots & \star \\
	& 1 & & (0) \\
	& & \ddots & \\
	&(0)  & & 1
\end{smallmatrix}\right)\bigg\}. \]

Denote by $\chi \coloneq \kappa(f)$ the linear form on $\g$ associated to $f$ and by $\overline{\chi} \coloneqq \chi|_{\g_{\geqslant 1}}$ its restriction to $\g_{\geqslant 1}$. Because this Lie algebra is Abelian, $\overline{\chi}$ is a character. The restriction map $\g^\vee \twoheadrightarrow (\g_{\geqslant 1})^\vee$ is a moment map, and we can twist it by the character $\overline{\chi}$, leading to the moment map
\[ \mu \colon \g^\vee \longrightarrow (\g_{\geqslant 1})^\vee, \quad \xi \longmapsto \xi|_{\g_{\geqslant 1}} + \overline{\chi}. \]
Its fibre over $0$ is the following affine subspace of $\g^\vee$:
\[ \mu^{-1}(0) = - \chi + (\g^\vee)_{\geqslant 0}. \]

To compute its quotient modulo the action of $G_{\geqslant 1}$, one introduces the Slodowy slice associated to the minimal $\sl_2$-triple:
\[ S_f \coloneqq - \chi + [\g, e]^\perp = \kappa(-f + \g^e). \]
The Slodowy slice is a closed subvariety of $ \mu^{-1}(0)$ because $\g^e$ is non-negatively graded.

\begin{theorem}[{\cite[Lemma 2.1]{gan2002quantization}}] \label{theorem:gan-ginzburg-isomorphism}
	The action map
	\[ \phi \colon G_{\geqslant 1} \times S_f \longrightarrow - \chi + (\g^\vee)_{\geqslant 0}, \quad (g, \xi) \longmapsto g \cdot \xi  \]
	is a $G_{\geqslant 1}$-equivariant isomorphism of algebraic varieties. By taking quotient on both sides with respect to the $G_{\geqslant 1}$-actions, one gets an isomorphism
	\[ S_f \cong \mu^{-1}(0) \git G_{\geqslant 1}. \]
\end{theorem}

\subsection{Restriction to a distinguished open set}

Consider the distinguished open subscheme
\[ \g^\vee[e^{-1}] \coloneqq \Spec ((\Sym \g)[e^{-1}]). \]
It is a Poisson variety, stable by the action of $G_{\geqslant 1}$ because $e$ is $G_{\geqslant 1}$-invariant, and we have the induced moment map given by restriction of $\mu$:
\[ \mu_{e^{-1}} \colon \g^\vee[e^{-1}] \longrightarrow (\g_{\geqslant_1})^\vee. \]

\begin{proposition} \label{proposition:hamiltonian-reduction-open}
	The natural open inclusion 
	\[ (\mu_{e^{-1}})^{-1}(0) \hooklongrightarrow \mu^{-1}(0) \]
	descends to a Poisson isomorphism between the corresponding Hamiltonian reductions:
	\[ (\mu_{e^{-1}})^{-1}(0) \git G_{\geqslant 1} \cong \mu^{-1}(0) \git G_{\geqslant 1} \]
	
\end{proposition}

Concretely, this proposition relies on the easy observation that $\mu^{-1}(0)$ is contained in the distinguished open subscheme $\g^\vee[e^{-1}]$. Because later on we will deal with localisation in noncommutative settings, we give details for the proof in term of commutative algebra.

\begin{proof}
	Denote by $I_{\mu}$ the ideal of $\Sym \g$ that is generated by the elements $x_{1,2} + 1$ and $x_{1,i}$ for $3 \leqslant i \leqslant d$. Then:
	\[ (\mu_{e^{-1}})^{-1}(0) = \Spec ((\Sym\g) [e^{-1}] / I_{\mu} [e^{-1}]) \quad \text{and} \quad \mu^{-1}(0) = \Spec (\Sym\g / I_{\mu}). \]
	Because localisation is exact \cite[Proposition 2.5]{eisenbud1995commutative}, one has a natural algebra isomorphism
	\[ (\Sym\g) [e^{-1}] / I_{\mu} [e^{-1}] \cong ((\Sym\g) / I_{\mu}) [e^{-1}] \]

	The image of $e = x_{1,2}$ in the quotient $(\Sym\g) / I_{\mu}$ is equal to $-1$, so it is invertible. Therefore, the localisation homomorphism
	\[ (\Sym\g) / I_{\mu} \longrightarrow ((\Sym\g) / I_{\mu}) [e^{-1}] \]
	is an isomorphism by \cite[Exercise 2.9]{eisenbud1995commutative}. This shows the isomorphism
	\[ (\mu_{e^{-1}})^{-1}(0) \cong \mu^{-1}(0). \]
	
	This isomorphism is clearly $G_{\geqslant 1}$-equivariant, so we get isomorphic quotients. This isomorphism is Poisson because it is induced by the $G_{\geqslant 1}$-equivariant Poisson embedding $\g^\vee[e^{-1}] \hookrightarrow \g^\vee$ that is compatible with the moment maps $\mu_{e^{-1}}$ and $\mu$.
\end{proof}

\subsection{Slodowy slice as the quotient of a distinguished open set} \label{subsection:slodowy-quotient-open}

The Lie algebra $\Z$-grading
\[ \g = \g_{-1} \oplus \g_0 \oplus \g_1 \]
is equivalent to an action of the multiplicative group $\Gm$ on $\g$ by Lie algebra automorphisms. Concretely, we have an algebraic group homomorphism
\[ \tau \colon \Gm \longrightarrow \operatorname{PGL}_d \cong G, \quad t \longmapsto \left(\begin{smallmatrix}
	t & & & 0 \\
	& 1 & & \\
	& & \ddots & \\
	0 & & & 1
\end{smallmatrix}\right) \in \mathrm{PGL}_d, \]
such that for any $i$ in $\{-1, 0, 1\}$, $x$ in $\g_i$ and $t$ in $\Gm$, one has:
\[ \tau(t) \cdot x = t^i x. \]

Define the following Abelian group of matrices:
\[ K \coloneqq \Gm \times \Ga{}^{d-2}. \]
We have the following closed embedding of algebraic group:
\[ \alpha \colon K \longrightarrow \mathrm{PGL}_d \cong G, \quad (t, a_1, \dots, a_{d-2}) \longmapsto \left(\begin{smallmatrix}
	t & & & & (0) \\
	0 & 1 & &  & \\
	0 & a_1 & 1  & & \\
	\vdots & \vdots & & \ddots &  \\
	0 & a_{d-2} & (0) & & 1
\end{smallmatrix}\right). \]
The corresponding Lie algebra $\kfk$ is spanned in $\g$ by $H$ and the $x_{i,2}$ for $3 \leqslant i \leqslant d$.

\begin{lemma} \label{lemma:afk}
	There is a vector space decomposition
	\[ \g_0 = \kfk \oplus (\g_0)^f. \]	
\end{lemma}

\begin{proof}
	A straightforward computation on the generators of $\kfk$ shows that $\kfk \cap (\g_0)^f = 0$. Because the grading is good in the sense of \cite{elashvili2005classification,brundan2007good}, there is an exact sequence
	\[ 0 \longrightarrow (\g_0)^f \longrightarrow \g_0 \overset{[f, \bullet]}{\longrightarrow} \g_{-1} \longrightarrow 0, \]
	therefore $\dim \g_0 = \dim(\g_0)^f + \dim \g_1$. Since $\dim \g_1 = \dim \kfk = d - 1$, we can conclude that the inclusion $\kfk \oplus (\g_0)^f \subseteq \g_0$ is an equality.
\end{proof}

For $1 \leqslant i \leqslant d - 2$ an integer, introduce the linear form $\chi_i \coloneqq \kappa(x_{i+2, 1})$. Then we have the following explicit description for the distinguished open set:
\[ \g^\vee[e^{-1}] = (\C - \{0\}) \chi + \sum_{i=1}^{d-2} \C \chi_i + (\g^\vee)_{\geqslant 0}. \]

\begin{lemma} \label{lemma:dilation-isomorphism}
	The action map
	\[ \begin{split}	
		\phi_0 \colon K \times (- \chi + (\g^\vee)_{\geqslant 0}) & \longrightarrow  \g^\vee[e^{-1}], \\
		((t, a_1, \dots, a_{d-2}), \xi) & \longmapsto  \alpha(t, a_1, \dots, a_{d-2}) \cdot \xi. 
	\end{split} \]
	is an isomorphism of varieties.
\end{lemma}

\begin{proof}
	The action of $(t, a_1, \dots, a_{d-2})$ in $K$ on $f$ is given by the formula
	\[ \alpha(t, a_1, \dots, a_{d-2}) \cdot (-f) = \left(\begin{smallmatrix}
		0 & & & & (0) \\
		-1/t & 0 & &  & \\
		-a_1 / t & 0 & 0  & & \\
		\vdots & \vdots & & \ddots &  \\
		-a_{d-2} / t & 0 & (0) & & 0
	\end{smallmatrix}\right). \]
	By applying $\kappa$, we get an isomorphism
	\[ K \overset{\sim}{\longrightarrow} (\g_{\geqslant 1})^\vee[e^{-1}], \quad (t, a_1, \dots, a_{d-2}) \longmapsto - \alpha(t, a_1, \dots, a_{d-2}) \cdot \overline{\chi}. \]
	Because the projection $\g^\vee \twoheadrightarrow (\g_{\geqslant 1})^\vee$ is $K$-equivariant, the pull-back of this isomorphism is the isomorphism stated in the lemma.
\end{proof}

The previous lemma implies that the adjoint action of the group $K \ltimes G_{\geqslant 1}$ stabilises the open subset $\g^\vee[e^{-1}]$. It acts by Poisson automorphisms, so the categorical quotient
\[ \g^\vee[e^{-1}] \git (K \ltimes G_{\geqslant 1}) \coloneq \Spec \big((\Sym \g)[e^{-1}]\big)^{K \ltimes G_{\geqslant 1}} \] 
inherits a Poisson structure from $\g^\vee$.

The following statement is an analogue of the Gan--Ginzburg isomorphism (\Cref{theorem:gan-ginzburg-isomorphism}), but the codomain is the open subscheme $\g^\vee[e^{-1}]$ instead of the moment map fibre $\mu^{-1}(0)$.

\begin{proposition} \label{proposition:dilation-gan-ginzburg}
	There is a $K \ltimes G_{\geqslant 1}$-equivariant isomorphism of varieties given by:
	\[ \begin{split}
		\phi' \colon (K \ltimes G_{\geqslant 1}) \times S_f & \longrightarrow \g^\vee[e^{-1}], \\
		((t, a_1, \dots, a_{d-2}), g, \xi) & \longmapsto \alpha(t, a_1, \dots, a_{d-2}) g \cdot \xi. 
	\end{split} \]
\end{proposition}

\begin{proof}
	Recall the isomorphism given by \Cref{theorem:gan-ginzburg-isomorphism}:
	\[ \phi \colon G_{\geqslant 1} \times S_f \overset{\sim}{\longrightarrow} - \chi + (\g^\vee)_{\geqslant 0}, \quad (g, \xi) \longmapsto g \cdot \xi.  \]
	Moreover, we have the isomorphism $\phi_0$ given by \Cref{lemma:dilation-isomorphism}. Therefore, we have the explicit formula for an inverse of $\phi'$:
	\[ \phi'^{-1}(\nu) = ((t, a_1, \dots, a_{d-2}), g, \xi) \quad \text{where} \quad \begin{cases}
		(t, a_1, \dots, a_{d-2}, \nu_1) \coloneq \phi_0{}^{-1}(\nu), \\
		(g, \xi) \coloneq \phi^{-1}(\nu_1).
	\end{cases} \]
\end{proof}

We can deduce from the previous proposition a description of the Slodowy slice as the categorical quotient of the distinguished open set.

\begin{theorem} \label{theorem:slodowy-git}
	There is an isomorphism
	\[ \mu^{-1}(0) \git G_{\geqslant 1} \cong \g^\vee[e^{-1}] \git (K \ltimes G_{\geqslant 1}) \]
	of Poisson varieties.
\end{theorem}

\begin{proof}
	Consider the inclusions
	\[ S_f \hooklongrightarrow \mu^{-1}(0) \hooklongrightarrow \g^\vee[e^{-1}] . \]
	After taking quotients, we get homomorphisms
	\[ S_f \longrightarrow \mu^{-1}(0) \git G_{\geqslant 1} \longrightarrow \g^\vee[e^{-1}]\git (K \ltimes G_{\geqslant 1}) \]
	of schemes that are isomorphisms because $\phi$ and $\phi'$ are isomorphisms. 
	
	To check they are Poisson, let us examine their comorphisms. By \Cref{proposition:hamiltonian-reduction-open}, there is a Poisson isomorphism
	\[ (\Sym \g / I_{\mu})^{G_{\geqslant 1}} \cong ((\Sym \g)[e^{-1}] / I_{\mu} [e^{-1}])^{G_{\geqslant 1}}   \]
	induced by localisation.
	
	Moreover, the natural algebra homomorphism
	\[ (\Sym \g)[e^{-1}] \longrightarrow  (\Sym \g)[e^{-1}] / I_{\mu} [e^{-1}] \]
	induces an algebra isomorphism
	\[ \big((\Sym \g)[e^{-1}]\big)^{K \ltimes G_{\geqslant 1}} \cong \big((\Sym \g)[e^{-1}] / I_{\mu} [e^{-1}]\big)^{G_{\geqslant 1}} \]
	after descending to invariant subalgebras. This isomorphism is Poisson because both Poisson brackets are induced by the one on $(\Sym \g)[e^{-1}]$.
\end{proof}

%\begin{corollary}
%	There is an isomorphism
%	\[ \Jet_\infty \mu^{-1}(0) \git \Jet_\infty G_{\geqslant 1} \cong \Jet_\infty \g^\vee[e^{-1}] \git \Jet_\infty(K \ltimes G_{\geqslant 1}) \]
%	of Poisson vertex schemes.
%\end{corollary}
%
%\cite[Theorem 5.1]{bozec2025functorial}
%
%Consider the nilpotent Lie subalgebra
%

\subsection{Decomposition theorem for the minimal Slodowy slice}

Consider the Poisson schemes
\begin{multline*}
	(\kfk \ltimes \g_{\geqslant 1})^\vee[e^{-1}] = \Spec (\Sym(\kfk \ltimes \g_{\geqslant 1})[e^{-1}]) \\ \text{and} \quad S_f \cong \g^\vee[e^{-1}] \git (K \ltimes G_{\geqslant 1}) = \Spec\big((\Sym \g)[e^{-1}]\big)^{K \ltimes G_{\geqslant 1}}. 
\end{multline*}  
Denote by $\overline{\chi}'$ the restriction of the linear form $\chi = (f \mid \bullet)$ to the Lie subalgebra $\kfk \ltimes \g_{\geqslant 1}$.

\begin{lemma} \label{lemma:open-orbit}
	The orbit map of $-\overline{\chi}'$
	\[ K \ltimes G_{\geqslant 1} \longrightarrow (\kfk \ltimes \g_{\geqslant 1})^\vee, \quad g \longmapsto - g \cdot \overline{\chi}' \]
	induces a scheme isomorphism
	\[ K \ltimes G_{\geqslant 1} \cong (\kfk \ltimes \g_{\geqslant 1})^\vee[e^{-1}]. \]
	Therefore $(\kfk \ltimes \g_{\geqslant 1})^\vee[e^{-1}]$ is an open coadjoint orbit, in particular a symplectic smooth variety.
\end{lemma}

\begin{proof}
	We compute the orbit map in stages, starting with the action of $G_{\geqslant 1}$. Because the group is unipotent, we can use the exponential of an element $x$ in $\g_{\geqslant 1}$. If we iterate the coadjoint action of $x$ on $\overline{\chi}'$, then $x \cdot (x \cdot \overline{\chi}') = 0$ because $[x, [x,f]]$ belongs to $\g_1$ and is orthogonal to $\kfk \ltimes \g_{\geqslant 1}$. 
	
	Moreover, because $(\bullet \mid \bullet)$ gives a perfect pairing between $\g_{-1}$ and $\g_1$, and by using the decomposition of \Cref{lemma:afk}, we deduce that 
	\[ (\bullet \mid \bullet) \colon \g_1 \otimes \g_0 \longrightarrow \C, \quad x \otimes y \longmapsto ([x, f] \mid y) = (x \mid [f, y]) \]
	induces a perfect pairing between $[\g_{\geqslant 1}, f]$ and $\kfk$. 
	
	Therefore, the $\g_1 \cong G_{\geqslant 1}$-orbit map induces an isomorphism:
	\[ \g_1 \overset{\sim}{\longrightarrow} -\overline{\chi}' + \kfk^\vee \subseteq (\kfk \ltimes \g_{\geqslant 1})^\vee, \quad x \longmapsto -\overline{\chi}' - x \cdot \overline{\chi}',  \]
	where $\kfk^\vee$ is identified to the subspace $(\g_{\geqslant 1})^\perp$ of $(\kfk \ltimes \g_{\geqslant 1})^\vee$.
	
	If we now look at the action of $K$, we can deduce the lemma from the computation of \Cref{lemma:dilation-isomorphism}.
%	We use the isomorphism $\kappa\colon \g \overset{\sim}{\rightarrow} \g^\vee$ to make the computation.
%	
%	The action of $(t, a_1, \dots, a_{n-2})$ in $A$ on $f$ is explicitly described by the formula
%	\[ \left(\begin{smallmatrix}
%		t & b_1 & b_2& \cdots & b_{d-1} \\
%		0 & 1 & &  & \\
%		0 & a_1 & 0  & & \\
%		\vdots & \vdots & & \ddots &  \\
%		0 & a_{d-2} & (0) & & 1
%	\end{smallmatrix}\right) \cdot f = \left(\begin{smallmatrix}
%		0 & & & & (0) \\
%		1/t & 0 & &  & \\
%		a_1 / t & 0 & 0  & & \\
%		\vdots & \vdots & & \ddots &  \\
%		a_{n-2} / t & 0 & (0) & & 0
%	\end{smallmatrix}\right). \]
%	By applying $\kappa$, we get the lemma.
\end{proof}

The restriction map $\g^\vee \twoheadrightarrow (\kfk \ltimes \g_{\geqslant 1})^\vee$ gives rise to a Poisson homomorphism
\[ \g^\vee[e^{-1}] \longrightarrow (\kfk \ltimes \g_{\geqslant 1})^\vee[e^{-1}], \quad \xi \longmapsto \xi|_{\kfk \ltimes \g_{\geqslant 1}}. \]
Consider the fibred product of this restriction map and the quotient map
\[ \g^\vee[e^{-1}] \longrightarrow \g^\vee[e^{-1}] \git (K \ltimes G_{\geqslant 1}), \quad \xi \longmapsto [\xi]_{\bmod K \ltimes G_{\geqslant 1}}, \]
that we denote by
\[ \begin{aligned}
	\delta \colon \g^\vee[e^{-1}] & \longrightarrow (\kfk \ltimes \g_{\geqslant 1})^\vee[e^{-1}] \times \g^\vee[e^{-1}] \git (K \ltimes G_{\geqslant 1}) \\
	\xi & \longmapsto  (\xi|_{\kfk \ltimes \g_{\geqslant 1}}, [\xi]_{\bmod K \ltimes G_{\geqslant 1}})
\end{aligned} \]

\begin{lemma} \label{lemma:gan-ginzburg-decomposition}
	The homomorphism of schemes
	\[ \begin{aligned}
		\psi \colon (K \ltimes G_{\geqslant 1}) \times S_f & \longrightarrow (\kfk \ltimes \g_{\geqslant 1})^\vee[e^{-1}] \times \g^\vee[e^{-1}] \git (K \ltimes G_{\geqslant 1}) \\
		(g,\xi) & \longmapsto  (g \cdot \xi|_{\kfk \ltimes \g_{\geqslant 1}}, [g \cdot \xi]_{\bmod K \ltimes G_{\geqslant 1}})
	\end{aligned} \]
	is an isomorphism.
\end{lemma}

\begin{proof}
	Because the grading is good, we have the inclusions $\kfk \ltimes \g_{\geqslant 1} \subseteq \g_{\geqslant 0} \subseteq [\g, e]$. Then, we get the inclusions $[\g, e]^\perp \subseteq (\g_{\geqslant 0})^{\perp} \subseteq (\kfk \ltimes \g_{\geqslant 1})^{\perp}$.	Therefore, for any $\xi$ in $S_f = - \chi + [\g, e]^\perp$ and for any $g$ in $K \ltimes G_{\geqslant 1}$, we have the equality:
	\[ g \cdot \xi|_{\kfk \ltimes \g_{\geqslant 1}} = - g \cdot \overline{\chi}'. \]
	Therefore, we can simplify the expression for $\psi$:
	\[ \psi(g, \xi) = (- g \cdot \overline{\chi}', [\xi]_{\bmod K \ltimes G_{\geqslant 1}}). \]
	By \Cref{lemma:open-orbit} and \Cref{theorem:slodowy-git}, we see that $\psi$ is an isomorphism because it is the product of the isomorphisms
	\[ K \ltimes G_{\geqslant 1} \cong (\kfk \ltimes \g_{\geqslant 1})^\vee[e^{-1}] \quad \text{and} \quad S_f \cong \g^\vee[e^{-1}] \git (K \ltimes G_{\geqslant 1}). \]
\end{proof}

We can state and prove the decomposition theorem for the minimal Slodowy slice.

\begin{theorem} \label{theorem:decomposition-slodowy}
	The scheme homomorphism
	\[ 	\delta \colon \g^\vee[e^{-1}] \longrightarrow (\kfk \ltimes \g_{\geqslant 1})^\vee[e^{-1}] \times \g^\vee[e^{-1}] \git (K \ltimes G_{\geqslant 1}) \]
	is a Poisson isomorphism, where the Poisson structure on the right-hand side product is the product of Poisson structures.
\end{theorem}

\begin{proof}
	We start by proving that $\delta$ is compatible with the Poisson structures, and then we prove it is an isomorphism.
	
	Both Poisson algebras
	\[  (\Sym(\kfk \ltimes \g_{\geqslant 1}))[e^{-1}] \quad \text{and} \quad \big((\Sym \g)[e^{-1}]\big)^{K \ltimes G_{\geqslant 1}} \]
	are Poisson subalgebras of $(\Sym \g)[e^{-1}]$. They mutually Poisson-commute because of the following equalities of subalgebras in $(\Sym \g)[e^{-1}]$:
	\[ \begin{aligned}
		\big((\Sym \g)[e^{-1}]\big)^{K \ltimes G_{\geqslant 1}} & = \big((\Sym \g)[e^{-1}]\big)^{\kfk \ltimes \g_{\geqslant 1}} \\
		& = \big\{ F \in (\Sym \g)[e^{-1}] \quad | \quad \text{for all} ~  x \in \kfk \ltimes \g_{\geqslant 1},  ~   \{ x, F \} = 0 \big\} \\
		& = \big\{ F \in (\Sym \g)[e^{-1}] \quad | \quad  \{\Sym(\kfk \ltimes \g_{\geqslant 1})[e^{-1}], F \} = 0\big\}. 
	\end{aligned} \]
	\begin{itemize}
		\item The first equality is justified by the connectedness of the group $K \ltimes G_{\geqslant 1}$. 
		\item The second equality is because the Poisson bracket $\{ x, \bullet \}$ coincides with the infinitesimal action of $x$ in $\kfk \ltimes \g_{\geqslant 1}$.
		
		\item The third equality follows because Poisson commuting with $\Sym(\kfk \ltimes \g_{\geqslant 1})[e^{-1}]$ is equivalent to Poisson commuting with $\kfk \ltimes \g_{\geqslant 1}$, since this space generates the non-localised algebra and the Poisson bracket with $e^{-1}$ is given by the formula: $\{e^{-1}, \bullet\} = - e^{-2} \{e, \bullet\}$.
	\end{itemize}
	
	Therefore, the natural algebra homomorphism
	\[ \delta^* \colon \Sym(\kfk \ltimes \g_{\geqslant 1})[e^{-1}] \otimes \big((\Sym \g)[e^{-1}]\big)^{K \ltimes G_{\geqslant 1}} \longrightarrow (\Sym \g)[e^{-1}] \]
	given by tensoring the inclusion homomorphisms is Poisson. So, $\delta$ is a Poisson variety homomorphism. 
	
	To prove it is an isomorphism, we notice that we have the following commutative diagram:
	\[ \begin{tikzcd}[column sep = small]
		& (K \ltimes G_{\geqslant 1}) \times S_f \arrow[ld, "\phi'"', "\sim"'{sloped}] \arrow[rd, "\psi", "\sim"'{sloped}] & \\
		\g^\vee[e^{-1}] \arrow[rr, "\delta"'] & & (\kfk \ltimes \g_{\geqslant 1})^\vee[e^{-1}] \times \g^\vee[e^{-1}] \git (K \ltimes G_{\geqslant 1}).
	\end{tikzcd} \]
	Since both $\phi'$ and $\psi$ are isomorphisms (see \Cref{proposition:dilation-gan-ginzburg} and \Cref{lemma:gan-ginzburg-decomposition}), $\delta$ is an isomorphism by composition.
\end{proof}

\section{Filtrations, localisation and group actions} \label{section:localisation-filtration}

\subsection{Recollections on filtered vector spaces and algebras}

An \emph{increasing filtration} on a vector space $V$ is the data of subspaces $V_p$ of $V$, for all $p$ in $\Z$, such that $V_p \subseteq V_{p+1}$.

A filtration is called:
\begin{itemize}
	\item \emph{separated} if $0 = \bigcap_{p \in \Z} V_p$,
	
	\item \emph{bounded from below} if $V_p = 0$ for some $p$ in $\Z$,
	
	\item \emph{bounded from above} if $V_p = V$ for some $p$ in $\Z$,
	
	\item \emph{exhaustive} if $V = \bigcup_{p \in \Z} V_p$.
\end{itemize}

The associated graded vector space of a filtered vector space $V$ is 
\[ \gr V \coloneqq \bigoplus_{p \in \Z} V_p / V_{p - 1}. \]
Let $V$ and $W$ be two filtered vector spaces. A linear map $\phi \colon V \rightarrow W$ is called a \emph{filtered homomorphism} if $f(V_p) \subseteq W_p$ for all $p$ in $\Z$. Then we can define an associated graded homomorphism
\[ \gr \phi \colon \gr V \longrightarrow \gr W \]
in an obvious way.

The following fact is standard, see \cite[Lemma 6.8]{kamnitzer2022hamiltonian} for a proof.

\begin{lemma} \label{lemma:quantisation-iso}
	Let $V$ and $W$ be two filtered vector spaces such that
	\begin{enumerate}
		\item both filtrations are exhaustive
		
		\item the filtration on $V$ is separated,
		
		\item the filtration on $W$ is bounded from below.
	\end{enumerate} 
	Let $\phi \colon V \rightarrow W$ be a filtered homomorphism. Therefore, if $\gr \phi \colon \gr V \rightarrow \gr W$ is an isomorphism, then $\phi \colon V \rightarrow W$ is an isomorphism.
\end{lemma}

The assignment
\[ (V \otimes W)_p \coloneqq \sum_{i + j = p} V_i \otimes W_j, \quad \text{for} \quad p \in \Z, \]
defines a filtration on the tensor product $V \otimes W$. There is a natural isomorphism
\[ \gr V \otimes \gr W \cong \gr(V \otimes W). \]

\subsection{Localisation of a filtered algebra} \label{subsectin:localisation-filtered}

In this section, we recall how one can localise a positively filtered algebra according to \cite[Section 2.1]{kapranov1998noncommutative}. See also \cite[Section 6]{kamnitzer2022hamiltonian} or \cite{skoda2006noncommutative}.

Let $A$ be a unital associative algebra equipped with an increasing filtration by subspaces $A_p$ ($p$ is in $\Z$), such that $A_p = 0$ for $p < 0$, and we have the following multiplicative properties:
\begin{equation} \label{equation:multiplicative-properties} 1 \in A_0, \quad A_p \cdot A_q \subseteq A_{p + q} \quad \text{and} \quad [A_p, A_q] \subseteq A_{p + q - 1}, \end{equation}
for $p$ and $q$ in $\Z$. One sees that $A_0$ is a commutative subalgebra and the associated graded algebra
\[ \gr A \coloneqq \bigoplus_{p \in \Z} A_p / A_{p - 1} = \bigoplus_{p \geqslant 0} A_p / A_{p - 1} \]
is commutative. Moreover, $A_0$ identifies with the degree $0$ subalgebra of $\gr A$.  %Moreover, $A$ is equipped with two commuting left and right $A_0$-actions respectively by left and right multiplication. Let $e$ be in $A_0$ and denote by 
%\[ A[e^{-1}] \coloneqq A_0[e^{-1}] \otimes_{A_0} A \] 
%the left localisation with respect to $e$, that is a left $A_0[e^{-1}]$ and a right $A_0$-module.

The following lemma and its proof is inspired by \cite[Proposition 2.1.5]{kapranov1998noncommutative}

\begin{lemma} \label{lemma:ore}
	Any multiplicative subset $S$ in $A_0$ satisfies the left and right Ore condition in $A$:
	\begin{enumerate}
		\item \label{lemma:ore,item:left-frac} for any $a$ in $A$ and $s$ in $S$, there are $b$ in $A$ and $t$ in $S$ such that $t a = b s$,
		
		\item \label{lemma:ore,item:left-zero} if $a$ in $A$ and $s$ satisfy $a s = 0$, then there is some $t$ in $S$ such that $t a = 0$,
		
		\item for any $a$ in $A$ and $s$ in $S$, there are $b$ in $A$ and $t$ in $S$ such that $a t = s b$,
		
		\item if $a$ in $A$ and $s$ satisfy $s a = 0$, then there is some $t$ in $S$ such that $a t = 0$.
	\end{enumerate}
\end{lemma}

\begin{proof}
	We check only the left conditions \eqref{lemma:ore,item:left-frac} and \eqref{lemma:ore,item:left-zero}.
	
	Assume that $a$ is in $A_n$ for some nonnegative integer $n$ and $s$ is in $S$. Then for any positive integer $k$:
	\[ s^{k} a = \bigg(\sum_{i=0}^{k - 1} s^{k - 1 - i} \ad(s)^i a\bigg) s + \ad(s)^{k} a. \] 
	Note that $\ad(s)^{k} a$ belongs to $A_{n - k}$. So for $k = n + 1$,
	$\ad(s)^{n + 1} a = 0$ since $A_{-1} = 0$, and we get
	\[ s^{n+1} a = \bigg(\sum_{i=0}^{n} s^{n - i} \ad(s)^i a\bigg) s. \]
	Then \eqref{lemma:ore,item:left-frac} holds for
	\[ b \coloneqq \sum_{i=0}^{n} s^{n - i} \ad(s)^i a \quad \text{and} \quad t \coloneqq s^{n + 1}. \]
	
	Now assume that $a s = 0$. Then for any positive integer $k$:
	\[ s^k a = \ad(s)^k a. \]
	For $k = n + 1$, we get
	\[ s^{n + 1} a = \ad(s)^{n + 1} a = 0. \]
	Therefore \eqref{lemma:ore,item:left-zero} holds for $t = s^{n + 1}$.
\end{proof}

By the general theory of Ore localisation \cite[Chapter II]{stenstrom1975rings}, we get the following.

\begin{proposition} \label{proposition:fraction}
	For any multiplicative subset $S$ in $A_0$, there is an algebra $A[S^{-1}]$ equipped with an algebra homomorphism $\iota \colon A \rightarrow A[S^{-1}]$ satisfying the following universal property. For any algebra homomorphism $\psi \colon A \rightarrow B$ such that $\psi(S)$ lies in the group of units $B^\times$ of $B$, there is a unique algebra homomorphism $\widetilde{\psi} \colon A[S^{-1}] \rightarrow B$ such that $\widetilde{\psi} \circ \iota = \psi$.
	
	The algebra $A[S^{-1}]$ is a flat left and right $A$-module. The natural map
	\[ A_0[S^{-1}] \otimes_{A_0} A \longrightarrow A[S^{-1}] \]
	is an isomorphism of left $A_0[S^{-1}]$-module.
\end{proposition}

There is a natural filtration on $A[S^{-1}]$: because $A_0 \rightarrow A_0[S^{-1}]$ is flat, the chain of inclusions
\[ A_0 \subseteq A_1 \subseteq \cdots \subseteq A_p \subseteq \cdots \subseteq A \]
gives rise to a chain of inclusion
\[ A_0[S^{-1}] \subseteq A_0[S^{-1}] \otimes_{A_0} A_1  \subseteq \cdots \subseteq A_0[S^{-1}] \otimes_{A_0} A_p \subseteq \cdots \subseteq A[S^{-1}] \]
after tensoring by $A_0[S^{-1}]$ and using the identification $A_0[S^{-1}] \otimes_{A_0} A \cong A[S^{-1}]$. Therefore, the assignment
\[ (A[S^{-1}])_p \coloneqq A_0[S^{-1}] \otimes_{A_0} A_p, \quad \text{for} \quad p \in \Z,  \]
defines an increasing filtration on $A[S^{-1}]$ by vector subspaces. This filtration clearly satisfies the multiplicative properties \eqref{equation:multiplicative-properties}. In particular, $\gr(A[S^{-1}])$ is a graded commutative algebra and contains $A_0[S^{-1}]$ as the degree 0 subalgebra.

\begin{proposition} \label{proposition:gr-localisation}
	There is a natural algebra isomorphism
	\[ (\gr A)[S^{-1}] \longrightarrow \gr(A[S^{-1}]). \]
\end{proposition}

\begin{proof}
	The natural homomorphism $A \rightarrow A[S^{-1}]$ is compatible with the filtrations on both sides, therefore we get an algebra homomorphism
	\[ \gr A \longrightarrow \gr(A[S^{-1}]). \]
	Under this last homomorphism, $S$ on the left-hand side is mapped to the degree $0$ component of the right-hand side, that is to say $A_0[S^{-1}]$. Therefore, by universal property we get the natural algebra homomorphism
	\begin{equation} \label{equation:graded-iso} (\gr A)[S^{-1}] \longrightarrow \gr(A[S^{-1}]). \end{equation}
	
	To show that it is an isomorphism, we note that, for any $p$ in $\Z$, the exact sequence
	\[ 0 \longrightarrow A_{p - 1} \longrightarrow A_p \longrightarrow \gr_p A \longrightarrow 0 \]
	gives rise to an exact sequence
	\[ 0 \longrightarrow (A[S^{-1}])_{p - 1} \longrightarrow (A[S^{-1}])_p \longrightarrow A_0[S^{-1}] \otimes_{A_0} (\gr_p A) \longrightarrow 0 \]
	after applying $A_0[S^{-1}] \otimes_{A_0} \bullet$. Therefore, we get an isomorphism
	\[ \gr_p(A[S^{-1}]) \cong A_0[S^{-1}] \otimes_{A_0} (\gr_p A) \]
	of $A_0[S^{-1}]$-modules. This isomorphism gives an inverse to the natural algebra homomorphism \eqref{equation:graded-iso}, that is then an isomorphism.
\end{proof}

%\begin{lemma}
%	The right multiplication by $e$ gives an invertible endomorphism of the left localisation $A[e^{-1}]$. Therefore, $A[e^{-1}]$ is naturally a left and right $A_0[e^{-1}]$-module.
%\end{lemma}
%
%
%\begin{proof}
%	The right multiplication by $e$ is given by the formula:
%	\[  x e = e x + [x, e], \quad \text{for} \quad x \in A. \]
%	It gives a left $A_0$-linear endomorphism of $A$, that extends to a $A_0[e^{-1}]$-linear automorphism of $A[e^{-1}]$. Moreover, $[\bullet, e]$ is nilpotent because the brackets goes down the filtration. As a sum of an invertible and a nilpotent endomorphism, the right multiplication by $e$ induces an invertible endomorphism of $A[e^{-1}]$. This proves the statement.
%\end{proof}
%
%\begin{proposition}
%	There is a unique algebra structure on $A[e^{-1}]$ such that the natural homomorphism
%	\[ A \longrightarrow A[e^{-1}] \]
%	is an algebra homomorphism.
%\end{proposition}
%
%\begin{proof}
%	By localisation, the multiplication on $A$ gives rise to a map
%	
%	
%	It is clearly bilinear and associative.
%\end{proof}

\subsection{Filtered representations of an algebraic group}

\label{subsection:filtered-representation}

In this section, we recall from \cite{juillard2025phd} the notion of a filtered module over a graded algebraic group, and some useful properties to study the subspace of invariants. See also \cite[Section 2.4]{hamilton2026relative} and \cite[Section 3.1]{hoskins2025moduli}, where very similar notions are introduced.

Let $N$ be an algebraic group and $\n$ be its Lie algebra. The data of a Lie algebra $\Z$-grading on $\n$ is equivalent to the data of an algebraic action of $\Gm$ on $\n$ by Lie algebra automorphisms. By analogy, we define a \emph{grading on $N$} as the data of an action of $\Gm$ on $N$ by algebraic group automorphisms.

If $X$ is an affine scheme equipped with an action $\sigma \colon \Gm \times X \rightarrow X$, then $\C[X]$ is graded by
$$ \C[X](p) = \{F \in \C[X] ~ | ~  \sigma^*(F) = t^p \otimes F \} ~ \quad ~ \text{for} \quad p \in \Z. $$
Therefore, if $N$ is a graded group, $\C[N]$ is a graded algebra.

Let $N$ be a graded group and $V$ be an $N$-module, that is to say there is a co-action $V \rightarrow \C[N] \otimes V$ \cite[Chapter 4, Section a]{milne2017algebraic}. We say that an $N$-module $V$ is a \emph{graded $N$-module} if $V = \bigoplus_{p \in \Z} V(p)$ is a graded vector space and if the co-action $V \rightarrow \C[N] \otimes V$ is a graded homomorphism. If $N$ is connected, a graded $N$-module $V$ is equivalently defined as an $N$-module $V$ equipped with a $\Z$-grading such that, for any $i$ and $p$ in $\Z$, $\n_i \cdot V(p) \subseteq V(i + p)$.

By analogy, we call a \emph{filtered $N$-module} a filtered vector space $V$ equipped with a co-action $V \rightarrow \C[N] \otimes V$ that is a filtered homomorphism. If $N$ is connected, a filtered $N$-module $V$ is equivalently defined as an $N$-module $V$ equipped with a filtration such that, for any $i$ and $p$ in $\Z$, $\n_i \cdot V_p \subseteq V_{i + p}$. The filtration on $V$ induces a filtration on $V^N$ by setting $(V^N)_p \coloneq V_p \cap V^N$.

\begin{lemma} \label{lemma:gr-action}
	If $N$ is graded and $V$ is a filtered $N$-module, then $\gr V$ has an induced graded $N$-module structure.
\end{lemma}

\begin{proof}
	Taking the graded homomorphism associated with the filtered co-action gives the graded $N$-module structure on $\gr V$. 
\end{proof}

Assume that $N$ is connected. Recall that for an $N$-module $V$, the level $0$ of the Lie algebra cohomology computes the $N$-invariants:
$$ \Hgy^0(\n, V) = V^\n = V^N. $$
This cohomology is defined as the cohomology of the Chevalley--Eilenberg complex
\[ \Cpx^\bullet(\n, V) \coloneq V \otimes \operatorname{Alt}^\bullet \n^\vee. \] 
See \cite[Section 7.7]{weibel1994introduction} for precise definitions.

The following proposition is a generalisation of \cite[Proposition 5.2]{gan2002quantization} and \cite[Proposition 4.1.2.5]{juillard2025phd}. It gives a sufficient condition for the operations $\gr$ and $(\bullet)^N$ for a filtered module over a unipotent group $N$.

\begin{proposition} \label{proposition:gr-invariant}
	Let $N$ be a graded unipotent group whose Lie algebra is nonpositively graded: $\n = \bigoplus_{i \leqslant 0} \n_i$. Let $V$ be a filtered $N$-module such that:
	\begin{enumerate}
		\item the filtration on $V$ is bounded from below and exhaustive,
		
		\item there is an isomorphism $\C[N] \otimes W \cong \gr V$ of $N$-modules, where $\C[N]$ is equipped with the left-multiplication action and $W$ denotes a trivial $N$-module.
	\end{enumerate}
	Then there is a natural graded linear isomorphism 
	$$ \gr (V^N) \overset{\sim}{\longrightarrow} (\gr V)^N. $$
	Moreover, the Lie algebra cohomology $\Hgy^n(\n, V)$ is zero whenever $n \neq 0$. 
\end{proposition}

\begin{remark}
	In \cite{juillard2025phd}, we assumed that the Lie algebra grading was negatively supported: $\n = \bigoplus_{i \leqslant -1} \n_i$. One can relax this hypothesis to a nonpositively supported grading, but the unipotence of the group is no longer automatic.
\end{remark}

\begin{proof}
	Consider the dual $\Z$-grading on $\n^\vee$, that is nonnegative: $\n^\vee = \bigoplus_{i \geqslant 0} (\n^\vee)_i$. For $p$ in $\Z$ and $n \geqslant 0$, define the linear subspace
	$$ \Cpx(\n, V)^n_p \coloneqq \sum_{\substack{\mathcal{l}, i_1, \dots, i_n \in \Z \\ \mathcal{l} + i_1 + \cdots + i_n \leqslant p}} V_{\mathcal{l}} \otimes ((\n^\vee)_{i_1} \wedge \cdots\wedge (\n^\vee)_{i_n}). $$ 
	This defines a filtration on the complex such that the differential satisfies the inclusion $d(C^\bullet_p) \subseteq C^\bullet_p$ for all $p$ in $\Z$. According to \cite[Section~5.4]{weibel1994introduction}, there is a spectral sequence~$\{E_{r}\}_{r=0}^\infty$ associated with this filtered cochain complex.	
	
	For $p, q$ in $\Z$, the $0$-page is given by:
	$$ E_0^{p, q} = \gr_p C^{q-p} \coloneqq \bigoplus_{\mathcal{l} + i_1 + \cdots + i_{q-p} = p} \gr_\mathcal{l} V \otimes ((\n^\vee)_{i_1} \wedge \cdots\wedge(\n^\vee)_{i_{q-p}}). $$
	It is the Chevalley--Eilenberg complex for the $\n$-module $\gr V$: for $n \in \Z$, there is the identification
	$$ \bigoplus_{p \in \Z}  E_0^{p, n+p} = \gr V \otimes \operatorname{Alt}^n \n^\vee. $$
	Hence the $1$-page is the Lie algebra cohomology of $\gr V$:
	$$  \bigoplus_{p \in \Z}  E_1^{p, n+p} = \Hgy^{n}(\n, \gr V). $$ 
	
	By using the $\n$-equivariant isomorphism $\C[N] \otimes W \cong \gr V$, we get
	$$ \Hgy^{n}(\n, \gr V) \cong W \otimes \Hgy^{n}(\n, \C[N]). $$
	But the cohomology $\Hgy^{n}(\n, \C[N])$ coincides with the algebraic de Rham cohomology of $N$ as a consequence of the isomorphism between the cotangent sheaf of $N$ and $\C[N] \otimes \n^\vee$ \cite[Proposition 4.4.2]{springer2009linear}. Because the exponential map $\exp \colon \n \overset{\sim}{\rightarrow} N$ provides an isomorphism between $N$ and the affine space $\n$, its cohomology is given by
	$$ \Hgy^{n}(\n, \C[N]) \cong \delta_{n = 0} \, \C. $$
	Hence, the first page $E_1^{p, q}$ is zero unless $p=q$: the spectral sequence collapses.
	
	The filtration on $C$ is exhaustive and bounded from below because so is the filtration on $V$ and the dual of the Lie algebra has only nonnegative degree components. Hence, according to \cite[Theorem 5.5.1]{weibel1994introduction}, the spectral sequence converges to the Lie algebra cohomology of $V$: the infinity page is given by
	$$ E_\infty^{p, q} = \gr_p \Hgy^{q-p}(C^\bullet). $$	
	Because of the collapsing we get isomorphisms $E_\infty^{p, q} \cong E_1^{p, q}$. In particular, taking $p=q$, we see that the natural embedding
	$$  \gr (V^N) = \gr (V^\n) \hooklongrightarrow (\gr V)^\n = (\gr V)^N$$
	induced by the inclusion $V^N \subseteq V$ is an isomorphism. Taking $q = p + n$ for nonzero $n$, we get the cohomology vanishing.
\end{proof}

The following proposition states that $\gr$ and $(\bullet)^{\Gm}$ always commute for a filtered module over the multiplicative group.

\begin{proposition} \label{proposition:weight-filtration}
	Let $V$ be a filtered $\Gm$-module, where $\Gm$ is trivially graded. Then $\gr V$ is equipped with a natural $\Gm$-module structure that preserves the grading. 
	
	Denote by $V_p(\delta)$ and $(\gr_p V)(\delta)$ the weight spaces of the $\Gm$-actions, where $\delta$ and $p$ are in $\Z$. There is a linear isomorphism
	$$ \gr_p(V(\delta)) \cong (\gr_p V)(\delta) $$
	induced by the inclusion $V(\delta) \subseteq V$.
	
	For $\delta = 0$, this gives the isomorphism
	\[ \gr_p(V^\Gm) \cong (\gr_p V)^{\Gm}.  \]
\end{proposition}

\begin{proof}
	The first part of the statement is a consequence of \Cref{lemma:gr-action} applied to $\Gm$, seen as a trivially graded group so that $V$ is a filtered $\Gm$-module. 
	
	The filtration on $V(\delta)$ is defined by
	$$ V(\delta)_p \coloneqq V(\delta) \cap V_p = V_p(\delta). $$
	The inclusion $V(\delta) \subseteq V$ respects the filtrations so we get a graded inclusion~$\gr (V(\delta)) \subseteq (\gr V) (\delta)$. But $\gr_p V$ decomposes in two possible ways:
	$$ \bigoplus_{\delta \in \Z} \gr_p (V(\delta)) = \gr_p V = \bigoplus_{\delta \in \Z}  (\gr_p V) (\delta), $$
	hence the  inclusion $\gr_p (V(\delta)) \subseteq (\gr_p V) (\delta)$ is an equality. 
\end{proof}

\section{Decomposition theorem in the quantum setting}
\label{section:quantum-setting}

\subsection{Finite W-algebra for the minimal orbit}

Finite W-algebras are defined by a quantum Hamiltonian reduction process \cite{kostant1978whittaker,premet2002special,gan2002quantization,brundan2007good}. See \cite[Section 7.1]{gan2006almost} for generalities on quantum Hamiltonian reduction.

Let $\U(\g)$ be the universal enveloping algebra of $\g$. It is equipped with an adjoint action of $G$ by algebra automorphisms and an induced adjoint action of $\g$. The map
\[ \widetilde{\mu} \colon \g_{\geqslant 1} \longrightarrow \U(\g), \quad a \longmapsto a + \chi(a) 1 \]
is a \emph{quantum comoment map} for the $G_{\geqslant 1}$-action, meaning that it is a Lie algebra homomorphism such that the following equality holds:
\[ [\widetilde{\mu}(a), F] = a \cdot F, \quad \text{for} \quad a \in \g_{\geqslant 1} \quad \text{and} \quad F \in \U(\g), \]
where $a \cdot F$ denotes the adjoint action.

The \emph{minimal finite W-algebra} is the quantum Hamiltonian reduction
\[ U(\g, f) \coloneqq (\U(\g) / I_{\widetilde{\mu}})^{G_{\geqslant 1}}, \]
where $I_{\widetilde \mu}$ is the left ideal of $\U(\g)$ that is generated by the elements $x_{1,2} + 1$ and $x_{1,i}$ for $3 \leqslant i \leqslant d$, that is to say by the image of $\widetilde{\mu}$. By a standard fact on quantum Hamiltonian reduction, $\U(\g, f)$ inherits an associative algebra structure from $\U(\g)$.

\subsection{The Kazhdan filtration on the W-algebra} We recall the definition of the Kazhdan filtration on $\U(\g)$ and $\U(\g, f)$ from \cite[Section 4]{gan2002quantization}. The semisimple element $H$ acts on $\g$ by derivations and its action is diagonalisable, defining the good grading. Therefore, we get a well-defined diagonalisable action on $\U(\g)$. The \emph{Kazhdan filtration} is defined by
\[ \U(\g)_p \coloneqq \sum_{n - m \leqslant p} \operatorname{Span} \{ x_1 \cdots x_n \in \U(\g) \quad | \quad [H, x_1 \cdots x_n] = m \, x_1 \cdots x_n \}, \]
where $p$ belongs to $\Z$.

Because $\g = \g_{-1} \oplus \g_0 \oplus \g_1$ and $e$ is in $\g_1$, $e$ is in the level $0$ for this filtration and $\U(\g)_p$ is non-zero only if $p \geqslant 0$. 

\begin{remark} \label{remark:nonnegative-filtration}
	The fact that the filtration is supported on $\Z_{\geqslant 0}$ comes from the fact that the good grading on $\g$ is supported on $\{-1, 0, 1\}$. If $f$ is not a minimal nilpotent element, it is almost always the case that this Kazhdan grading contains components with negative degrees. For example, take the regular nilpotent element in $\sl_3$.
\end{remark}

The Kazhdan filtration induces filtrations on the quotient $\U(\g) / I_{\widetilde{\mu}}$ and on the finite W-algebra $\U(\g, f) = (\U(\g) / I_{\widetilde{\mu}})^{G_{\geqslant 1}}$. They are both called \emph{Kazhdan filtrations} too. The following properties are easy to prove, see \cite{gan2002quantization} or \cite[Section 4.2.1]{juillard2025phd}.

\begin{lemma}
	\begin{enumerate}
		\item The Kazhdan filtration is exhaustive.
		
		\item It satisfies the multiplicative conditions \eqref{equation:multiplicative-properties}, so $\gr \U(\g)$ is commutative.
		
		\item There is an isomorphism $\C[\g^\vee] = \Sym \g \cong \gr \U(\g)$. 
		
		\item The commutative subalgebra $\U(\g)_0$ is a free commutative algebra generated by the subspace $\g_1$:
		\[ \U(\g)_0 \cong \Sym \g_1. \]
	\end{enumerate}
\end{lemma} 

Because $e$ belongs to $\U(\g)_0$, we can apply the localisation procedure described in \Cref{subsectin:localisation-filtered}. Then, $\U(\g)[e^{-1}]$ is well-defined, it inherits a Kazhdan filtration and there is a Poisson algebra isomorphism
\[ \C[\g^\vee[e^{-1}]] = (\Sym \g)[e^{-1}] \cong \gr \big(\U(\g)[e^{-1}]\big). \]

The following theorem proves that the Slodowy slice is quantised by the finite W-algebra.

\begin{theorem}[{\cite[Theorem 4.1]{gan2002quantization}}] \label{theorem:slodowy-quantisation}
	The isomorphism $\Sym \g \cong \gr \U(\g)$ induces algebra isomorphisms
	\[ \Sym \g / I_\mu \cong \gr(\U(\g) / I_{\widetilde{\mu}}) \quad \text{and} \quad \C[S_f] \cong (\Sym \g / I_\mu)^{G_{\geqslant 1}} \cong \gr \U(\g, f), \]
	and the second isomorphism is Poisson.
\end{theorem}

\subsection{The W-algebra is a subalgebra of the localisation}

We state a quantisation of \Cref{theorem:slodowy-git}, that is the reconstruction of the finite W-algebra as a subalgebra of $\U(\g)[e^{-1}]$. Because the punctured line $\C^\times e$ is stable under the adjoint action of $K \ltimes G_{\geqslant 1}$, the action of this group extends from $\U(\g)$ to the localisation $\U(\g)[e^{-1}]$.

If we consider the opposite grading $\g = \bigoplus_{i \in \Z}\g_{-i}$, then the Lie subalgebra $\kfk \ltimes \g_{\geqslant 1}$ is $\Z_{\leqslant 0}$-graded and $\U(\g)[e^{-1}]$ is a filtered module because the Kazhdan filtration satisfies the property:
\[ [x, \big(\U(\g)[e^{-1}]\big)_p] \subseteq \big(\U(\g)[e^{-1}]\big)_{1 - i + p - 1} = \big(\U(\g)[e^{-1}]\big)_{p - i}, \]
for $x$ in $\g_i \subseteq \U(\g)_{1 - i}$ and $p$ an integer.

\begin{theorem} \label{theorem:finite-w-algebra-coset}
	There is an algebra isomorphism
	\[ \big(\U(\g)[e^{-1}]\big)^{K \ltimes G_{\geqslant 1}} \cong \U(\g, f).  \]
\end{theorem}

To prove this theorem, we start by computing the associated graded algebra of $\big(\U(\g)[e^{-1}]\big)^{K \ltimes G_{\geqslant 1}}$ for the induced Kazhdan filtration.

\begin{proposition} \label{proposition:gr-invariant-localisation}
	The natural algebra homomorphism
	\[ \gr \big(\big(\U(\g)[e^{-1}]\big)^{K \ltimes G_{\geqslant 1}}\big) \longrightarrow \big(\gr \big(\U(\g)[e^{-1}]\big)\big)^{K \ltimes G_{\geqslant 1}} \]
	is an isomorphism. Therefore, we have a Poisson algebra isomorphism
	\[ \gr \big(\big(\U(\g)[e^{-1}]\big)^{K \ltimes G_{\geqslant 1}}\big) \cong \C[\g^\vee[e^{-1}] \git (K \ltimes G_{\geqslant 1})]. \]
\end{proposition}

\begin{proof}
	The group $K = \Gm \times \Ga{}^{d-2}$ decomposes as the product of a torus and a unipotent group.
	
	To prove the proposition, we can first deal with the action of the unipotent normal subgroup $\Ga{}^{d-2} \ltimes G_{\geqslant 1}$ in $K \ltimes G_{\geqslant 1}$. To show the isomorphism
	\begin{equation} \label{equation:gr-invariant-unipotent} \gr \big(\big(\U(\g)[e^{-1}]\big)^{\Ga{}^{d-2} \ltimes G_{\geqslant 1}}\big) \cong \big(\gr \big(\U(\g)[e^{-1}]\big)\big)^{\Ga{}^{d-2}\ltimes G_{\geqslant 1}}, \end{equation}
	we want to use \Cref{proposition:gr-invariant}. 
	
	According to \Cref{proposition:dilation-gan-ginzburg}, we have the $\Ga{}^{d-2} \ltimes G_{\geqslant 1}$-equivariant isomorphism
	\[ \C[\Ga{}^{d-2} \ltimes G_{\geqslant 1}] \otimes \C[\Gm \times S_f] \cong (\Sym \g)[e^{-1}] \cong \gr (\U(\g)[e^{-1}]).  \]
	With the discussion on the Lie algebra grading at the beginning of the section, we have all the assumptions to apply \Cref{proposition:gr-invariant}, and we get the isomorphism~\eqref{equation:gr-invariant-unipotent}.
	
	The isomorphism \eqref{equation:gr-invariant-unipotent} is $\Gm$-equivariant, so we get the isomorphism
	\[ \big(\gr \big(\big(\U(\g)[e^{-1}]\big)^{\Ga{}^{d-2} \ltimes G_{\geqslant 1}}\big)\big)^{\Gm} \cong \big(\gr \big(\U(\g)[e^{-1}]\big)\big)^{K \ltimes G_{\geqslant 1}}. \]
	between the $\Gm$-invariant subalgebras. By applying \Cref{proposition:weight-filtration} on the left-hand side, we get the desired isomorphism:
	\[ \gr \big(\big(\U(\g)[e^{-1}]\big)^{K \ltimes G_{\geqslant 1}}\big) \cong \big(\gr \big(\U(\g)[e^{-1}]\big)\big)^{K \ltimes G_{\geqslant 1}}. \]
	
	The second isomorphism of the proposition follows because we have the Poisson isomorphism
	\[ \gr \big(\U(\g)[e^{-1}]\big) \cong (\Sym \g)[e^{-1}] = \C[\g^\vee[e^{-1}]] \]
	which is $K \ltimes G_{\geqslant 1}$-equivariant. By taking the invariant subalgebras, we get:
	\[ \gr \big(\big(\U(\g)[e^{-1}]\big)^{K \ltimes G_{\geqslant 1}}\big) \cong \big(\gr \big(\U(\g)[e^{-1}]\big)\big)^{K \ltimes G_{\geqslant 1}} \cong \C[\g^\vee[e^{-1}] \git (K \ltimes G_{\geqslant 1})]. \]
\end{proof}

Now, we can prove the main theorem.

\begin{proof}[Proof of \Cref{theorem:finite-w-algebra-coset}]
	First, the localisation homomorphism $\U(\g) \hookrightarrow \U(\g)[e^{-1}]$ induces an algebra isomorphism
	\[ (\U(\g) / I_{\widetilde \mu} )^{G_{\geqslant 1}} \cong (\U(\g)[e^{-1}] / I_{\widetilde \mu} [e^{-1}])^{G_{\geqslant 1}}. \]
	The proof is similar to that of \Cref{proposition:hamiltonian-reduction-open}: we use the exactness of localisation and the fact that $e$ is already invertible in the quotient $\U(\g) / I_{\widetilde \mu}$, so we have an isomorphism
	\[ \U(\g) / I_{\widetilde \mu}  \cong \U(\g)[e^{-1}] / I_{\widetilde \mu} [e^{-1}] \]
	according to \cite[Proposition 3.6]{stenstrom1975rings}.
	
	Second, there is a natural algebra homomorphism
	\[ \lambda\colon (\U(\g)[e^{-1}])^{K \ltimes G_{\geqslant 1}} \longrightarrow  (\U(\g)[e^{-1}] / I_{\widetilde \mu} [e^{-1}])^{G_{\geqslant 1}} \cong \U(\g, f). \]
	Both sides are equipped with Kazhdan filtrations that are bounded from below and exhaustive. So if $\gr \lambda$ is an isomorphism, then \Cref{lemma:quantisation-iso} implies that $\lambda$ is an isomorphism.
	
	By \Cref{proposition:gr-invariant-localisation} and \Cref{theorem:slodowy-quantisation}, we have the algebra isomorphisms
	\[ \begin{multlined}
		\gr \big(\big(\U(\g)[e^{-1}]\big)^{K \ltimes G_{\geqslant 1}}\big) \cong \C[\g^\vee[e^{-1}] \git (K \ltimes G_{\geqslant 1})] \\ \text{and} \quad \gr \U(\g, f) \cong \C[\mu^{-1}(0) \git G_{\geqslant 1}].
	\end{multlined} \]
	Under these isomorphisms, $\gr \lambda$ is conjugate to the isomorphism
	\[ \C[\g^\vee[e^{-1}] \git (K \ltimes G_{\geqslant 1})] \cong  \C[\mu^{-1}(0) \git G_{\geqslant 1}] \]
	given by \Cref{theorem:slodowy-git}. Therefore, $\gr \lambda$, and then $\lambda$, are isomorphisms.	
\end{proof}

\subsection{Decomposition theorem for finite W-algebras}

\begin{lemma}
	The two subalgebras
	\[ \U(\kfk \ltimes \g_{\geqslant 1})[e^{-1}] \quad \text{and} \quad \big(\U(\g)[e^{-1}]\big)^{K \ltimes G_{\geqslant 1}} \cong \U(\g, f)  \]
	of $\U(\g)[e^{-1}]$ commute with each other.
\end{lemma}

\begin{proof}
	The lemma follows immediately from the following equalities:
	\[ \begin{aligned}
		\big(\U(\g)[e^{-1}]\big)^{K \ltimes G_{\geqslant 1}} & = \big(\U(\g)[e^{-1}]\big)^{\kfk \ltimes \g_{\geqslant 1}} \\
		& = \big\{ F \in \U(\g)[e^{-1}] \quad | \quad \text{for all} ~  x \in \kfk \ltimes \g_{\geqslant 1},  ~   [x, F] = 0 \big\} \\
		& = \big\{ F \in \U(\g)[e^{-1}] \quad | \quad  [\U(\kfk \ltimes \g_{\geqslant 1})[e^{-1}], F]  = 0\big\}. 
	\end{aligned} \]
	\begin{itemize}
		\item The first equality is justified by the connectedness of the group $K \ltimes G_{\geqslant 1}$. 
		\item The second equality is because the bracket $[x, \bullet]$ coincides with the infinitesimal action of $x$ in $\kfk \ltimes \g_{\geqslant 1}$.
		
		\item The third equality follows because commuting with $\U(\kfk \ltimes \g_{\geqslant 1})[e^{-1}]$ is equivalent to commuting with $\kfk \ltimes \g_{\geqslant 1}$, since this space generates the non-localised algebra and the bracket with $e^{-1}$ is given by the formula: $[e^{-1}, \bullet] = - e^{-1} [e, \bullet] e^{-1}$.
	\end{itemize}
\end{proof}

Consequently, there is a natural algebra homomorphism
\[ \widetilde{\delta} \colon \U(\kfk \ltimes \g_{\geqslant 1})[e^{-1}] \otimes \big(\U(\g)[e^{-1}]\big)^{K \ltimes G_{\geqslant 1}} \longrightarrow \U(\g)[e^{-1}] \]
given by multiplication. This provides the decomposition isomorphism for the finite W-algebra.

\begin{theorem} \label{theorem:decomposition-quantum}
	The algebra homomorphism
	\[ \widetilde{\delta} \colon \U(\kfk \ltimes \g_{\geqslant 1})[e^{-1}] \otimes \big(\U(\g)[e^{-1}]\big)^{K \ltimes G_{\geqslant 1}} \longrightarrow \U(\g)[e^{-1}] \]
	is an isomorphism. Therefore, we have an algebra isomorphism
	\[ \U(\kfk \ltimes \g_{\geqslant 1})[e^{-1}] \otimes \U(\g, f) \cong \U(\g)[e^{-1}]. \]
\end{theorem}

\begin{proof}
	Both sides are equipped with induced filtrations coming from the Kazhdan filtration on $\U(\g)$ and localising at $e$ in $\U(\kfk \ltimes \g_{\geqslant 1})_0 \subseteq \U(\g)_0$. Because the filtrations on both sides are bounded from below and exhaustive, we just need to prove that
	\[ \gr \widetilde\delta \colon \gr \big(\U(\kfk \ltimes \g_{\geqslant 1})[e^{-1}]\big) \otimes \gr \big(\big(\U(\g)[e^{-1}]\big)^{K \ltimes G_{\geqslant 1}}\big) \longrightarrow \gr\big( \U(\g)[e^{-1}]\big) \]
	is an isomorphism to apply \Cref{lemma:quantisation-iso}.
	
	By applying \Cref{proposition:gr-localisation}, we obtain the algebra isomorphisms:
	\[ \gr\big(U(\kfk \ltimes \g_{\geqslant 1}\big)[e^{-1}]) \cong \Sym(\kfk \ltimes \g_{\geqslant 1})[e^{-1}]  \]
	and 
	\[ \gr\big(\U(\g)[e^{-1}]\big) \cong \Sym(\g)[e^{-1}]. \]
	
	Moreover, \Cref{proposition:gr-invariant-localisation} gives the algebra isomorphism
	\[ \gr \big(\big(\U(\g)[e^{-1}]\big)^{K \ltimes G_{\geqslant 1}}\big) \cong \big((\Sym\g)[e^{-1}]\big)^{K \ltimes G_{\geqslant 1}}.  \]
	
	Under these isomorphisms, $\gr \widetilde{\delta}$ is conjugate to the Poisson algebra isomorphism 
	\[ \delta^* \colon \Sym(\kfk \ltimes \g_{\geqslant 1})[e^{-1}] \otimes \big((\Sym \g)[e^{-1}]\big)^{K \ltimes G_{\geqslant 1}}  \overset{\sim}{\longrightarrow} (\Sym \g)[e^{-1}]  \]
	provided by \Cref{theorem:decomposition-slodowy}. This concludes the proof.
\end{proof}

\section{Vertex algebras, filtration, localisation and Lie algebras action} \label{section:vertex-algebra}

\subsection{Reminders on graded vertex algebras} \label{subsection:vertex-algebras}

For the details about the general theory of vertex algebras, see~\cite{kac1998vertex,frenkel2004vertex,desole2006finite}. We recall some basic notations and facts here.

Let $V$ be a vertex superalgebra over the complex number field $\C$. Denote by $1$ the vacuum vector, which is a non-zero even vector in $V$. Denote by $\partial : V \rightarrow V$ the translation operator, which is an even operator. The state-field correspondence is an even $\C$-linear homomorphism which maps $a \otimes b$ in $V \otimes V$ to a Laurent series $a(z) b$ in~$V(\!(z)\!)$. It is denoted by $a(z)b \coloneqq \sum_{n \in \Z} a_{(n)}b\, z^{-n-1}$ and the coefficient $a_{(n)}b$ is called the $n$-th product of $a$ and $b$. The formal power series
\[ a(z) \coloneqq \sum_{n \in \Z} a_{(n)} \quad \text{in }\quad (\operatorname{End} V)[\![z, z^{-1}]\!] \]
is called the \emph{field} associated to $a$.

The fundamental relation satisfied by the these fields is the Borcherds identity, see \cite[Section 3.3.10]{frenkel2004vertex} or \cite[(1.28)]{desole2006finite}:
\begin{equation} \label{equation:borcherds}
	\begin{multlined}
	\sum_{j \geqslant 0} (-1)^j \binom{n}{j} \big(a_{(m + n - j)}(b_{(k + j)} c) - (-1)^n (-1)^{|a||b|} b_{(n + k - j)} (a_{(m + j) c} \big) \\
	= \sum_{j \geqslant 0} \binom{m}{j} (a_{(n + j)} b)_{(m + k - j)} c,
	\end{multlined}  
\end{equation}
for any parity-homogeneous vectors $a, b, c$ in $V$ and any $m, n, k$ in $\Z$. We use the notation $|a|$ for the parity of $a$, that is to say $0$ if it is even, $1$ if it is odd.

We recall that the \emph{normally ordered product} is defined as $\NO{a b} \coloneqq a_{(-1)} b$. In particular, the vacuum satisfies $\NO{a 1} = \NO{a 1} = a$, and for any nonnegative integer $n$:
\[ \NO{\partial^n(a) b} = n! \, a_{(-n - 1)} b. \]
The normally ordered is non-associative in general, so we adopt the convention
\[ \NO{a b c} \coloneqq \NO{a \NO{b c}}, \]
where $a$, $b$ and $c$ are in $V$.

The \emph{operator product expansion (OPE)} of two elements $a$ and $b$ is the data of products $a_{(n)} b$ for all integers $n \geqslant 0$. It is denoted by:
\[ a(z_1) b(z_2) \sim \sum_{n \geqslant 0} \frac{1}{(z_1 - z_2)^n} (a_{(n)} b)(z_2). \]
Two element $a$ and $b$ are said to \emph{commute with each} other if and only if \[  [a(z_1), b(z_2)] = 0 \quad \text{in} \quad (\operatorname{End} V)[\![z_1, (z_1)^{-1}, z_2, (z_2)^{-1}]\!]. \]
By the Borcherds identity, it is equivalent to:
\[ a_{(n)} b = 0\quad \text{for all} \quad n \geqslant 0. \]

%
%For $a$ and $b$ in $V$, the {normally ordered product} is $\NO{a b} \coloneqq a_{(-1)} b$. For~$a, b, c$ in $\V$, we adopt the convention: $ \NO{a b c} \coloneqq  \NO{a (\NO{b c})}$. For $n$ a nonnegative integer, there is the relation $a_{(-n-1)} b = \frac{1}{n!} \NO{\partial(a) b}$. The $\lambda$-bracket is the formal sum defined by~$[a_\lambda b] \coloneqq \sum_{n=0}^\infty (a_{(n)}b)\frac{\lambda^n}{n!}$ in $V \otimes_\C \C[\lambda]$, it controls the nonassociativity and the noncommutativity of the normally ordered product.

In this paper, we call a \emph{graded vertex superalgebra} a vertex superalgebra $V$ such that there is a super vector space decomposition
\[ V = \bigoplus_{\Delta \in \Z_{\geqslant 0}} V(\Delta) \]
satisfying the inclusions
\[ V(\Delta){}_{(n)} V(\Delta') \subseteq V(\Delta + \Delta' - n - 1) \]
for any $n$ in $\Z$, and $\Delta, \Delta'$ in $\Z_{\geqslant 0}$. In particular $1$ must belong to $V(0)$ because $1_{(-1)}1 = 1$.

\begin{remark}
	It makes sense to consider a vertex algebra grading which is not supported on $\Z_{\geqslant 0}$. In this paper, such graded vertex algebras will not appear.
\end{remark}

\begin{theorem}[{\cite[Theorem 4.14]{li2004vertex}}] \label{theorem:li-filtration}
	Let $V$ be a graded vertex superalgebra. The data of subspaces
	\[ V_p \coloneqq \operatorname{Span} \Bigg\{ a_1{}_{(-n_1 - 1)} \cdots a_r{}_{(-n_r - 1)} 1 \quad \Bigg| \quad \begin{cases}
		r \in \Z_{\geqslant 0}, \\
		n_1, \dots, n_r \in \Z_{\geqslant 0}, \\
		\Delta_1, \dots, \Delta_r \in \Z_{\geqslant 0}, \\
		\Delta_1 + \cdots + \Delta_r \leqslant p, \\
		a_\bullet \in V(\Delta_\bullet)
	\end{cases} \Bigg\} \]
	for $p$ in $\Z_{\geqslant 0}$ defines an increasing, lower bounded and exhaustive filtration. For any integers $n$, $p$ and $q$, we have the embeddings:
	\begin{equation} \label{equation:li-filtration-property} \begin{aligned}
		V_p{}_{(n)} V_q & \subseteq V_{p + q} \quad \text{for} \quad n < 0, \\
		V_p{}_{(n)} V_q & \subseteq V_{p + q - 1} \quad \text{for} \quad n \geqslant 0.
	\end{aligned} \end{equation}
	Therefore, the associated graded vector space $\gr V$ inherits a super commutative vertex superalgebra structure from $V$ and a \emph{vertex Poisson structure}.
\end{theorem}

As in \cite{arakawa2026faithful}, we call this filtration the \emph{increasing Li filtration}.

\begin{remark}
	In \cite[Definition 2.7 and Theorem 2.12]{li2005abelianizing}, Li introduces an upper bounded decreasing filtration $V_{\mathrm{Li}}^\bullet$, called the \emph{decreasing Li filtration} in \cite{arakawa2026faithful}, such that $\gr_{\mathrm{Li}} V$ is a vertex Poisson superalgebra. In fact, there is an isomorphism
	\[ \gr_{\mathrm{Li}} V \cong \gr V \]
	of vertex Poisson superalgebras \cite[Proposition 2.6.1]{arakawa2012remark}.
\end{remark}

%\begin{proposition}[{\cite[Proposition 2.6.1]{arakawa2012remark}}] \label{proposition:compare-filtrations}
%	Let $V$ be a graded vertex superalgebra. Then, for any integer $p$, $V_p$ and $V_{\mathrm{Li}}^p$ are graded subspaces and there is the equality
%	\[ V_{\mathrm{Li}}^p (\Delta) = V_{\Delta - p} (\Delta). \]
%	This equalities induce 
%\end{proposition} 

Let $\{a_i\}_{i=1}^n$ be a set of odd or even vectors in $V$. Denote by $\leqslant$ the lexicographic order on $\{1, \ldots, n\} \times \Z_{\geqslant 0}$. The vertex superalgebra $V$ is said \emph{freely generated} by the \emph{strong generators} $\{a_i\}_{i=1}^n$ if the products
\begin{equation*}
	a_{i_1}{}_{(-n_1 - 1)} \cdots a_{i_r}{}_{(-n_r - 1)} 1 \quad \text{where} \quad 
	\begin{cases}
		& r \geqslant 0 \quad \text{is an integer}, \\
		& i_\bullet \in \{1, \ldots, n\} \quad \text{and} \quad n_\bullet \in \Z_{\geqslant 0}, \\
		& (i_j, n_j) \leqslant (i_{j+1}, n_{j+1}) \quad \text{if} \quad a_{i_j} \quad \text{is even}, \\
		& (i_j, n_j) < (i_{j+1}, n_{j+1}) \quad \text{if} \quad a_{i_j} \quad \text{is odd},
	\end{cases}
\end{equation*}
form a basis of $V$ as a vector space. 

\subsection{Arc schemes and commutative vertex algebras} 
\label{subsection:arc-schemes}

Standard references for the theory of jet and arc schemes are \cite{ein2009jet,chambert2018motivic}. In the context of vertex algebras, see \cite[Section 2.1]{arakawa2015associated}.

Let $R$ be a (super)commutative unital associative superalgebra that is of finite type. Denote by $\J_\infty R$ the universal differential superalgebra generated by $R$, it comes with a natural superalgebra embedding $R \hookrightarrow \J_\infty R$ and satisfies the universal property given by \cite[Remark 3.1]{ein2009jet}.

\begin{lemma}[{\cite[Section 1.4]{frenkel2004vertex}}]
	Like any supercommutative differential algebra, $\J_\infty R$ is a supercommutative vertex algebra.
\end{lemma} 

Let $X \coloneqq \Spec R$ be the associated affine super scheme. The affine super scheme  $\J_\infty X \coloneqq \Spec \J_\infty R$ represents the functor 
\[ \text{superalgebra} ~ A \longmapsto \{\text{superalgebra homomorphisms} ~  R \rightarrow  A[\![t]\!]\}, \]
where $A[\![t]\!]$ is the superalgebra of formal power series with coefficients in $A$. The scheme $\J_\infty X$ is called the \emph{arc scheme} of $X$.

\begin{proposition}[{\cite[Proposition 2.3]{arakawa2012remark}}]
	If $R$ is a Poisson superalgebra, the Poisson bracket extends to a unique vertex Poisson super bracket on $\J_\infty R$.
\end{proposition}

For any nonnegative integer $m$, the functor
\[ \text{superalgebra} ~ A \longmapsto \{\text{superalgebra homomorphisms} ~  R \rightarrow  A[t]/(t^{m+1})\} \]
is represented by an affine scheme $\J_m X = \Spec \J_m R$, where $A[t]/(t^{m+1})$ denotes the algebra $A \otimes \C[t] / t^{m + 1} \C[t]$. The scheme $\J_m X$ is called the \emph{$m$-jet scheme} of $X$. They come with natural super scheme homomorphisms
\[ \J_\infty X \longrightarrow \J_{m+1} X \longrightarrow \J_{m} X, \]
called the truncation morphisms, making $\J_\infty X$ the projective limit of the $\J_m X$, for $m$ in $\Z_{\geqslant 0}$.

The jet and arc constructions commute with localisation in the following sense. If $S$ is a multiplicative subset of $R = \C[X]$, then there are natural algebra isomorphisms:
\[ \J_m(R[S^{-1}]) \cong (\J_mR)[S^{-1}] \quad \text{and} \quad \J_\infty(R[S^{-1}]) \cong (\J_\infty R)[S^{-1}].  \]

By construction, if $X$ is of finite type then $\J_m X$ is always of finite type for any finite $m$. But it is not the case for $\J_\infty X$ in general. If $X$ is a smooth irreducible variety then $\J_m X$ is a smooth irreducible variety for any finite $m$ \cite[Corollary~2.11]{chambert2018motivic}. 

If $N$ is an affine group scheme acting on an affine scheme $X$, then for any finite or infinite $m$, $\J_m N$ is an affine group scheme acting on $\J_m X$ \cite[Remark 2.8]{ein2009jet}. 

\subsection{Lie algebra cohomology for a projective limit of Lie algebras}

In this section, we review the properties of modules of the projective limit of a system of algebraic groups. This is closely related to the notion of smooth module introduced in \cite[Section 1.9]{kazhdan1993tensor}.

Let $\{N_m\}_{m \in \Z_{\geqslant 0}}$ be a projective system of algebraic groups such that each homomorphism $N_{m + 1} \rightarrow N_m$ is \emph{dominant}, meaning that $\C[N_m] \rightarrow \C[N_{m+1}]$ is \emph{injective}. Denote by $N_\infty$ the affine group scheme given by the projective limit.

\begin{example}
	If $N$ is an algebraic group then for any finite $m$, the jet scheme $\J_m N$ is an algebraic group. The arc group $\J_\infty N$ is a group scheme \cite[Remark 2.8]{ein2009jet} and is the projective limit of the algebraic groups $\J_m N$. If $N$ is connected then $\J_m N$ is connected. 
	
	The Lie algebras of the jet and arc groups are given by
	\[ \operatorname{\mathfrak{Lie}} \J_m N = \J_m \n \coloneqq \n[t] / (t^{m+1}) \quad \text{and} \quad \operatorname{\mathfrak{Lie}} \J_\infty N = \J_\infty \n \coloneqq \n[\![t]\!], \]
	where the brackets are induced by extension of scalars (it is a well-known fact, see \cite[Lemma 5.1.6.1]{juillard2025phd} for a proof).
\end{example}

\begin{proposition} \label{proposition:colimit-module-pro-group}
	Let $V$ be an $N_\infty$-module, and denote by $\rho \colon V \rightarrow \C[N_\infty] \otimes V$ the corresponding co-action map. Then the subspaces
	\[ V_{(m)} \coloneqq \{v \in V \quad | \quad \rho(v) \in \C[N_m] \otimes V \}, \quad \text{where} \quad m \in \Z_{\geqslant 0}, \] 
	define an inductive system of submodules such that, for each $m$, the action of $N_\infty$ on $V_{(m)}$ factorises through the natural morphism $N_\infty \rightarrow N_m$. Moreover, $V$ is the colimit of the $V_{(m)}$.
	
	If $\phi \colon V \rightarrow W$ is an $N_\infty$-module homomorphism, then for all integers $m$, there is an inclusion
	\[ \phi(V_{(m)}) \subseteq W_{(m)}. \]
\end{proposition}

\begin{remark}
	In the previous statement, when we write $\rho(v) \in \C[N_m] \otimes V$, we identify $\C[N_m]$ with its image in $\C[N_\infty]$.
\end{remark}

\begin{proof}
	We start with the first statement. If $v$ belongs to $V_{(m)}$, we can write $\rho(v)$ as a finite sum
	\[ \rho(v) = \sum_i F_i \otimes v_i \]
	where the $F_i$ are functions in $\C[N_m]$ and the $v_i$ are in $V$. Using the compatibility of the comodule structure with the counit and comultiplications on $\C[N_m]$ and $\C[N_\infty]$ \cite[Formulae (24)]{milne2017algebraic}, we get that 
	\[ v = \sum_i v_i \quad \text{and} \quad \rho(v_i) \in \C[N_m] \otimes V. \]
	Therefore, we deduce that
	\[ \rho(V_{(m)}) \subseteq \C[N_m] \otimes V_{(m)}, \]
	so $V_{(m)}$ is an $N_\infty$-submodule such that the action of $N_\infty$ factorises through the natural homomorphism $N_\infty \rightarrow N_m$. 
	
	There are obvious inclusions $V_{(m)} \subseteq V_{(m+1)}$ by using the group homomorphisms $N_{m+1} \rightarrow N_m$. 
	Therefore we get a natural module homomorphism
	\[ \operatorname{colim}\limits_{m \geqslant 0} V_{(m)} \hooklongrightarrow V. \]
	Let us prove it is surjective. For any $v$ in $V$, we can write $\rho(v)$ as a finite sum
	\[ \rho(v) = \sum_i F_i \otimes v_i, \] 
	where the $F_i$ are functions in $\C[N_\infty]$ and the $v_i$ are in $V$. Because we have a filtered colimit, these finitely many $F_i$ must all lie in one $\C[N_m]$ for some $m$ that is big enough and then $v$ belongs to this $V_{(m)}$. This concludes the proof of the isomorphism
	\[ \operatorname{colim}\limits_{m \geqslant 0} V_{(m)} \cong V. \]
	
	Let us consider the second statement. We have the equation
	\[ (\id_{\C[N_\infty]} \otimes \phi) \circ \rho = \rho \circ \phi, \]
	by definition of a $\C[N_\infty]$-comodule homomorphism $\phi \colon V \rightarrow W$. If $v$ belongs to $V_{(m)}$ then 
	\[ \rho(\phi(v)) =  (\id_{\C[N_\infty]} \otimes \phi) \circ \rho(v) =  (\id_{\C[N_m]} \otimes \phi) \circ \rho(v) \]
	belongs to $\C[N_m] \otimes W$, so $\phi(v)$ is in $W_{(m)}$.
\end{proof}

Recall that the action by left multiplication of an affine group scheme $N$ on itself defines an $N$-module structure on the ring of functions $\C[N]$.

\begin{lemma} \label{lemma:filtration-function-progroup}
	For any nonnegative integer $m$, the image of the algebra embedding $\C[N_m] \hookrightarrow \C[N_\infty]$ equals $\C[N_\infty]_{(m)}$.
\end{lemma}

\begin{proof}
	Denote by $\mu$ the multiplication on the group scheme $N_\infty$. The comultiplication
	\[ \mu^*\colon \C[N_\infty] \longrightarrow \C[N_\infty] \otimes \C[N_\infty] \]
	is the co-action in this case. Denote by $\varepsilon \colon \C[N_\infty] \rightarrow \C$ the algebra homomorphism defined by evaluating polynomial functions at the unit of the group. By \cite[Proposition 3.6]{milne2017algebraic}:
	\[ (\id_{\C[N_\infty]} \otimes \varepsilon) \circ \mu^*(F) = F \otimes 1 \]
	for all $F$ in $\C[N_\infty]$. If we take $F$ in $\C[N_\infty]_{(m)}$, we deduce that $F$ belongs to $\C[N_m]$ by applying the previous identity.
	
	This gives the inclusion $\C[N_\infty]_{(m)} \subseteq \C[N_m]$. The other inclusion holds because $N_\infty \rightarrow N_m$ is a group homomorphism, so we have an equality.
\end{proof}

Denote by $\n_m$ the Lie algebra of $N_m$ and $\n_\infty$ the Lie algebra of $N_\infty$, which is the projective limit of the $\n_m$. For $V$ an $N_\infty$-module, \Cref{proposition:colimit-module-pro-group} gives submodules $V_{(m)}$ with an induced action of $\n_m$. 

As explained in \Cref{subsection:filtered-representation}, we define the Lie algebra cohomology $\Hgy^\bullet(\n_m, V_{(m)})$ as the cohomology of the Chevalley--Eilenberg complex. The cochain complex given by the colimit
\[ \Cpx^\bullet(\n_\infty, V) \coloneqq \operatorname{colim}\limits_{m \geqslant 0} \Cpx^\bullet(\n_m, V_{(m)}) = V \otimes \operatorname{Alt}^\bullet\big(\operatorname{colim}\limits_{m \geqslant 0} (\n_m)^\vee\big) \]
is the Chevalley--Eilenberg complex for $V$. The Lie algebra cohomology of $V$ is then defined as
\[ \Hgy^\bullet(\n_\infty, V) \coloneqq \Hgy(\Cpx^\bullet(\n_\infty, V)) = \operatorname{colim}\limits_{m \geqslant 0}\Hgy^\bullet(\n_m, V_{(m)}), \]
where the last equality comes from the fact that cohomology commutes with filtered colimits.

\begin{lemma} \label{lemma:coholology-pro-algebra}
	Let $V$ be an $N_\infty$-module. Its degree $0$ Lie algebra cohomology gives the Lie algebra invariants:
	\[ \Hgy^0(\n_\infty, V) = V^{\n_\infty}. \]

	If all the $N_m$ are connected, then 
	\[ V^{N_\infty} = V^{\n_\infty}. \]
\end{lemma}

\begin{proof}
	The first statement is given by the following equalities:
	\[ V^{\n_\infty} =  \operatorname{colim}\limits_{m \geqslant 0} (V_{(m)})^{\n_m} = \operatorname{colim}\limits_{m \geqslant 0} \Hgy^0(\n_m, V_{(m)}) = \Hgy^0(\n_\infty, V). \]
	
	Because each $N_m$ is connected and an algebraic group (in particular a connected smooth variety), one has the equality $(V_{(m)})^{N_m} = (V_{(m)})^{\n_m}$. Taking the colimits, the second part of the statement follows.
\end{proof}

%Let $V$ be a $\n_\infty$-module. By analogy with the terminology of \cite{kazhdan1993tensor}, we say that $V$ is a smooth module if $V$ can be written as the colimit of an inductive system of submodules $V_{(m)}$, where $m$ is a nonnegative integer, such that the action of $\n_\infty$ on $V$ factorises through through the natural morphism $\n_\infty \rightarrow N_m$.

\begin{proposition} \label{proposition:gr-invariant-pro}
	Assume that, for all nonnegative integers $m$, $N_m$ is a graded unipotent group, such that its Lie algebra is nonpositively graded $\n_m = \bigoplus_{i \leqslant 0} \n_{m,i}$. Moreover, assume that the homomorphism $N_{m+1} \rightarrow N_m$ is compatible with the gradings. Then, the projective limit $N_\infty$ inherits a grading.
	
	Let $V$ be an $N_\infty$-module, equipped with an increasing filtration, that is filtered in the sense of \Cref{subsection:filtered-representation}. Make the following assumptions:
	\begin{enumerate}
		\item the filtration on $V$ is bounded from below and exhaustive,

		\item there is an isomorphism $\C[N_\infty] \otimes W \cong \gr V$ of $N_\infty$-modules, where $\C[N_\infty]$ is equipped with the left-multiplication action and $W$ denotes a trivial $N_\infty$-module.
	\end{enumerate}
	Then there is a natural graded linear isomorphism 
	$$ \gr (V^{N_\infty}) \overset{\sim}{\longrightarrow} (\gr V)^{N_\infty}. $$
	Moreover, the Lie algebra cohomology $\Hgy^n(\n_\infty, V)$ is zero whenever $n \neq 0$. 
\end{proposition}

\begin{proof}
	The argument is similar to the proof of \Cref{proposition:gr-invariant}. The cohomology of $V$ is defined as the cohomology of the complex
	\[ \Cpx^\bullet(\n_\infty, V) = V \otimes \operatorname{Alt}^\bullet\big(\operatorname{colim}\limits_{m \geqslant 0} (\n_m)^\vee\big). \]
	The grading on the Lie algebra and the filtration on $V$ equip this cochain complex with an exhaustive filtration which is bounded from below. Therefore, there is a spectral sequence~$\{E_{r}\}_{r=0}^\infty$ associated with this filtered cochain complex \cite[Section~5.4]{weibel1994introduction}.
	
	The $1$-page of this complex is given by the Lie algebra cohomology of $\gr V$:
	\[ \bigoplus_{p \in \Z} E_1^{p, p + n} = \Hgy^{n}(\n_\infty, \gr V) = \operatorname{colim}\limits_{m \geqslant 0} \Hgy^{n}(\n_m, (\gr V)_{(m)}), \]
	where $n$ is a nonnegative integer.
	
	By using \Cref{proposition:colimit-module-pro-group} and \Cref{lemma:filtration-function-progroup}, for any $m$ we get an $N_m$-equivariant isomorphism 
	\[ \C[N_m] \otimes W \cong (\gr V)_{(m)}. \]
	Therefore, we can apply \Cref{proposition:gr-invariant} and we get that
	\[ \Hgy^{n}(\n_m, (\gr V)_{(m)}) = 0 \quad \text{for} \quad n \neq 0. \]
	By taking the colimit, we get
	\[ \Hgy^{n}(\n_\infty, (\gr V)) = 0 \quad \text{for} \quad n \neq 0, \]
	so the spectral sequence $\{E_{r}\}_{r=0}^\infty$ collapses at the $1$-page. Therefore we can conclude the argument in the same way as for \Cref{proposition:gr-invariant}.
\end{proof}

\begin{example} \label{example:grading-jet-groups}
	If $N$ is a graded group, then the jet groups $\{\J_m N\}_{m \geqslant 0}$ form a projective system of graded groups. If $\n = \bigoplus_{i \in \Z} \n_i$ denotes the induced Lie algebra grading, the Lie algebra grading on $\n[t]/(t^{m+1})$ is the grading induced by $\n$ by imposing that the formal variable $t$ is of degree $0$.
\end{example}

\subsection{Localisation \emph{à la} Borisov} \label{subsection:localisation-vertex}

We explain the process of localisation for a graded vertex algebra as it was originally developed in \cite{borisov2001vertex,malikov1999chiral}. Our main reference is \cite[Sections 2 and 3]{arakawa2026faithful}.

In this subsection, let $V$ be a $\Z_{\geqslant 0}$-graded vertex algebra. The normally ordered product equips $V(0)$ with a commutative associative algebra structure with unit given by the vacuum vector. The embedding $V(0) \hookrightarrow V$ extends uniquely to a vertex algebra homomorphism
\[ \J_\infty V(0) \longrightarrow V, \]
making $V$ a vertex algebra module over $\J_\infty V(0)$ (that is called a $V(0)$-loop module in \cite{borisov2001vertex}).

\begin{theorem}[{\cite[Propositions 2.3 and 3.5]{borisov2001vertex} and \cite[Theorem 3.12]{arakawa2026faithful}}] \label{theorem:borisov-localisation}
	For any multiplicative subset $S$ in $V(0)$, one can define a graded vertex algebra $V[S^{-1}]$ equipped with a graded vertex algebra homomorphism $\iota \colon V \rightarrow V[S^{-1}]$ satisfying the following universal property. 
	
	For any graded vertex algebra $W$ and any graded vertex algebra homomorphism $\psi \colon V \rightarrow W$ such that $\psi(S)$ lies in the group of units $W(0)^\times$ of $W(0)$, there is a unique graded vertex algebra homomorphism $\widetilde{\psi} \colon V[S^{-1}] \rightarrow W$ such that $\widetilde{\psi} \circ \iota = \psi$.	
\end{theorem}

The vertex algebra $V[S^{-1}]$ is called the localisation of $V$ with respect to the multiplicative subset $S$. Such a localisation can also be defined for any vertex algebra module $M$ over $V$ that is (non-negatively) graded. One denotes by $\iota_M \colon M \rightarrow M[S^{-1}]$ this localisation. It satisfies an analogous universal property as $V[S^{-1}]$, see \cite[Proposition 2.3]{borisov2001vertex} and \cite[Proposition 3.2]{arakawa2026faithful}. Denote by $V\text{\rm{-}}\mathcal{mod}$ the category of graded $V$-modules.

\begin{proposition}[{\cite[Lemma 3.3]{arakawa2026faithful}}]
	Localisation defines a functor
	\[ V\text{\rm{-}}\mathcal{mod} \longrightarrow V[S^{-1}]\text{\rm{-}}\mathcal{mod}, \quad M \longmapsto M[S^{-1}] \]
	which is an exact functor.
\end{proposition}

Because $V[S^{-1}]$ is graded, it is equipped with an increasing Li filtration and one can show that this filtration is the localisation of the Li filtration on $V$ \cite[Equation (3.10)]{arakawa2026faithful}.

\begin{proposition}[{\cite[Lemma 3.13]{arakawa2026faithful}}]  \label{proposition:gr-borisov-localisation}
	There is a natural vertex Poisson algebra isomorphism
	\[ (\gr V)[S^{-1}] \overset{\sim}{\longrightarrow} \gr(V[S^{-1}]), \]
	where the left-hand side denotes the usual localised commutative algebra.
\end{proposition}

\subsection{Universal affine vertex algebra and its localisation}

To the simple Lie algebra $\g$ and a complex number $k$, one can associate a vertex algebra $\V^k(\g)$, called the \emph{universal affine vertex algebra}. See \cite[Section 2.1]{frenkel1992vertex}, \cite[Section 2.4]{frenkel2004vertex} or \cite[Section 1.8]{desole2006finite}. The vertex algebra $\V^k(\g)$ contains $\g$ as a linear subspace. It is freely generated by any ordered basis in $\g$ \cite[Proposition 1.20]{desole2006finite}.

If $\g = \bigoplus_{i \in \Z} \g_i$ is graded, then we can equip $\V^k(\g)$ with a vertex algebra grading
\[ \V^k(\g) = \bigoplus_{\Delta \in \Z} \V^k(\g)(\Delta) \] such that $\g_i$ is included in $\V^k(\g)(1 - i)$ \cite[Remark 1.23]{desole2006finite}. In general, such a grading is not necessarily supported on $\Z_{\geqslant 0}$, but if we pick our choice of good grading $\g = \g_{-1} \oplus \g_0 \oplus \g_1$ for the minimal nilpotent $f$, then we get that
\begin{equation} \label{equation:grading-universal-affine} 
	\V^k(\g) = \bigoplus_{\Delta \in \Z_{\geqslant 0}} \V^k(\g)(\Delta), 
\end{equation}
and the degree $0$ commutative subalgebra is:
\[ \Sym \g_1 \cong \V^k(\g)(0). \]
This can be seen by picking an ordered basis of $\g$ that is homogeneous for the grading.

By \Cref{theorem:li-filtration}, $\V^k(\g)$ is equipped with the increasing Li filtration and $\gr \V^k(\g)$ is a commutative vertex algebra. By using the fact that $\V^k(\g)$ is freely generated by any basis of $\g$, one can show that the inclusion of $\g$ in $\V^k(\g)$ induces a vertex Poisson algebra isomorphism
\begin{equation} \label{equation:gr-universal-affine}
	\C[\J_\infty \g^\vee] = \J_\infty(\Sym \g) \cong \gr \V^k(\g).
\end{equation}

Because $e$ belongs to $\g_1$, we can use the Borisov localisation process given by \Cref{theorem:borisov-localisation}, giving a graded vertex algebra homomorphism
\[ \V^k(\g) \longrightarrow \V^k(\g)[e^{-1}]. \]
We deduce from \eqref{equation:gr-universal-affine} and \Cref{proposition:gr-borisov-localisation} the vertex Poisson algebra isomorphism
\[ \J_\infty(\Sym \g)[e^{-1}] = \J_\infty(\Sym \g[e^{-1}]) \cong \gr (\V^k(\g)[e^{-1}]). \]

\begin{remark}
	The localisation homomorphism $\V^k(\g) \rightarrow \V^k(\g)[e^{-1}]$ is injective because its associated graded map $\J_\infty(\Sym \g) \rightarrow \J_\infty(\Sym \g)[e^{-1}]$ is injective, and the increasing Li filtration is bounded from below.
\end{remark}

Any module $M$ over $\V^k(\g)$ carries an action of the arc algebra $\J_\infty \g = \g[\![t]\!]$, that is given by the formula
\[ x t^n \cdot a = x_{(n)} a, \]
where $x$ is in $\g$, $n$ is a nonnegative integer and $a$ is in $M$. If one takes the vertex subalgebra $\V^k(\mathfrak h)$ generated by a Lie subalgebra $\mathfrak{h}$ in $\g$, then we get an action of $\J_\infty \mathfrak h = \mathfrak h[\![t]\!]$.

\section{Chiral Hamiltonian reduction} \label{section:brst}

We recall the BRST formalism for chiral Hamiltonian reduction. Our main references are \cite{feigin1990quantization,kac2003quantum,arakawa2007representation,arakawa2015localization}. See also \cite[Section 3.5]{genra2025reduction} and references therein fore more details on the construction of the complex and the geometric interpretation.

\subsection{The fermionic vertex algebra}

Let $\V(\g_{\geqslant 1})$ be the vertex subalgebra of $\V^k(\g)$ that is generated by the subspace $\g_{\geqslant 1}$ of $\g$. It is a graded subalgebra for the grading \eqref{equation:grading-universal-affine}.

Denote by $\F(\g_{\geqslant 1} \oplus (\g_{\geqslant 1})^\vee)$ the \emph{fermionic vertex superalgebra} associated to the vector space $\g_{\geqslant 1} \oplus (\g_{\geqslant 1})^\vee$ equipped with the duality pairing \cite[Section 1.8]{desole2006finite}. It is a vertex algebraic analogue of a Clifford algebra, see \cite[Section 2.4.2]{genra2025reduction}. It contains $\g_{\geqslant 1} \oplus (\g_{\geqslant 1})^\vee$ as a linear subspace, and it is freely generated by any basis of $\g_{\geqslant 1} \oplus (\g_{\geqslant 1})^\vee$. If $a$ is an element of $\g_{\geqslant 1} \oplus (\g_{\geqslant 1})^\vee$, denote by $\phi(a)$ its image in the vertex superalgebra.

There is a vertex algebra grading
\[ \F(\g_{\geqslant 1} \oplus (\g_{\geqslant 1})^\vee) = \bigoplus_{\Delta \in \Z} \F(\g_{\geqslant 1} \oplus (\g_{\geqslant 1})^\vee)(\Delta) \] 
such that $\g_i$ is contained in $\F(\g_{\geqslant 1} \oplus (\g_{\geqslant 1})^\vee)(1 - i)$ and $(\g_i)^\vee$ is contained in $\F(\g_{\geqslant 1} \oplus (\g_{\geqslant 1})^\vee)(i)$ \cite[Equation (5.3)]{desole2006finite}. In the particular case of this good grading for minimal $f$, $\g_{\geqslant 1} = \g_1$ so we get a grading supported on $\Z_{\geqslant 0}$.

The vertex algebra $\F(\g_{\geqslant 1} \oplus (\g_{\geqslant 1})^\vee)$ can be equipped with the \emph{charge grading},
\[ \F(\g_{\geqslant 1} \oplus (\g_{\geqslant 1})^\vee) = \bigoplus_{n \in \Z} \F(\g_{\geqslant 1} \oplus (\g_{\geqslant 1})^\vee)^{\mathrm{ch} = n}, \] 
which is not a vertex algebra grading. Generators from $\g_{\geqslant 1}$ are of charge $-1$, generators from $(\g_{\geqslant 1})^\vee$ are of charge $+1$ and the translation operator preserves the charge.

\subsection{BRST cochain complex}

Let $\{a_i\}_{i=1}^{d-1}$ be a basis of $\g_{\geqslant 1}$ and $\{a_i^\vee\}_{i=1}^{d-1}$ be the dual basis for the duality pairing with $(\g_{\geqslant 1})^\vee$. The vertex superalgebra
\[ \V(\g_{\geqslant 1}) \otimes \F(\g_{\geqslant 1} \oplus (\g_{\geqslant 1})^\vee) \]
contains an element of charge $1$ and conformal degree $1$, called the \emph{BRST charge}:
\[ Q \coloneqq \sum_{i=1}^{d-1} a_i \otimes \phi(a_i^\vee) - \frac{1}{2} \sum_{i, j=1}^{d - 1} 1 \otimes \phi(a_i^\vee){}_{(-1)}\big(\phi(a_j^\vee){}_{(-1)}\phi([a_i, a_j])\big). \]
The zero-mode $Q_{(0)}$ is an endomorphism that squares to $0$, increases the charge-degree by $1$ and preserves the vertex algebra grading \cite[Sections 2.1 and 2.2]{kac2003quantum}. Then, the pair
\[ (\V(\g_{\geqslant 1}) \otimes \F(\g_{\geqslant 1} \oplus (\g_{\geqslant 1})^\vee), Q_{(0)}) \]
forms a cochain complex, called \emph{BRST cochain complex}, that is a graded vertex algebra. Each graded piece is a cochain subcomplex, so we have a vertex algebra grading.

Let $\widehat{\mu}_V \colon \V(\g_{\geqslant 1}) \rightarrow V$ be a (purely even) vertex algebra homomorphism, that is usually called a \emph{chiral comoment map}. We assume that $V$ is a graded vertex algebra and that $\widehat{\mu}_V$ is compatible with the gradings. The pair
\[ \big(\Cpx(V) \coloneqq V \otimes \F(\g_{\geqslant 1} \oplus (\g_{\geqslant 1})^\vee), \widehat{\mu}_V(Q)_{(0)}\big) \]
is a cochain complex that is a graded vertex algebra. Each graded piece is a cochain subcomplex. In the following proposition, denote by $[c]$ the cohomology class of some closed element $c$ in $\Cpx(V)$.

The following facts are standard, see for example \cite[Section 3.15]{arakawa2007representation}.

\begin{proposition} \label{proposition:brst-cohomology}
	The cohomology $\Hgy(\Cpx(V))$ enjoys the following properties:
	\begin{enumerate}
		\item it is a vertex superalgebra whose vacuum is $[1]$, translation operator is induced by the translation operator on the BRST complex and the state-field correspondence is given by:
		\[ [a](z)[b] \coloneqq [a(z)b] = \sum_{n \in \Z} [a_{(n)}b] \, z^{- n - 1} \quad \text{for} \quad [a], [b] \in \Hgy(\Cpx(V)),  \]
		
		\item it is a graded vertex algebra for the grading induced by the grading on the complex:
		\[ \Hgy(\Cpx(V)) = \bigoplus_{\Delta \in \Z_{\geqslant 0}}  \Hgy(\Cpx(V)(\Delta)), \]
		
		\item the degree $0$ cohomology $\Hgy^0(\Cpx(V))$ is a purely even vertex subalgebra.
	\end{enumerate}
\end{proposition}

This construction is functorial. Fix $\widehat{\mu}_V \colon \V(\g_{\geqslant 1}) \rightarrow V$ and $\widehat{\mu}_W \colon \V(\g_{\geqslant 1}) \rightarrow W$  two chiral comoment maps, and $\psi \colon V \rightarrow W$ a vertex algebra homomorphism such that $\widehat{\mu}_W = \psi \circ \widehat{\mu}_V$. Then there is a natural vertex super algebra homomorphism between the corresponding cochain complexes that descends to the corresponding cohomologies. 

\subsection{Li filtration on the quantum Hamiltonian reduction}

Let $V$ be a vertex algebra equipped with a chiral comoment map $\widehat{\mu}_V \colon \V(\g_{\geqslant 1}) \rightarrow V$. Let $X$ be an affine Poisson scheme with an action of $G_{\geqslant 1}$ and a moment map $\mu_X\colon X \rightarrow (\g_{\geqslant 1})^\vee$. We assume that $V$ is a \emph{chiral quantisation} of $X$ in the following sense: there is a vertex Poisson algebra isomorphism $\gr V \cong \C[\J_\infty X]$ that makes the triangle
\[ \begin{tikzcd}[column sep = small]
	\gr V \arrow[rr, "\sim"] &  & \C[\J_\infty X] \\
	& \arrow[lu, "\gr \widetilde{\mu}_V"] \arrow[ru, "(\J_\infty \mu_X)^*"'] \gr V(\g_{\geqslant 1})\cong \C[\J_\infty (\g_{\geqslant})^\vee]  &
\end{tikzcd} \]
commute.

Like any graded vertex superalgebra, the BRST complex $\Cpx(V)$ carries an increasing Li filtration. Because the differential preserves the vertex algebra grading and is a derivative for the modes of the vertex algebra, $\Cpx(V)$ is a filtered complex in the sense of \cite[Section 5.4]{weibel1994introduction} and $\gr \Cpx(V)$ has a natural structure of cochain complex.

There is an induced filtration on the cohomology vertex superalgebra $\Hgy(\Cpx(V))$ \cite[Section 5.5]{weibel1994introduction}: for $p$ in $\Z$, we denote by $\Hgy(\Cpx(V))_{(p)}$ the image of the homomorphism $\Hgy(\Cpx(V)_p) \rightarrow \Hgy(\Cpx(V))$. This defines an exhaustive filtration on $\Hgy(\Cpx(V))$ that is supported on $\Z_{\geqslant 0}$. 

A priori, this induced filtration is not the increasing Li filtration on the vertex superalgebra $\Hgy(\Cpx(V))$. However, this filtration inherits the properties \eqref{equation:li-filtration-property} from the Li filtration on the complex and we have the following lemma.

\begin{lemma}
	The associated graded vector space $\gr \Hgy(\Cpx(V))$ inherits a structure of vertex Poisson superalgebra.
\end{lemma}

Moreover, because both the induced filtration satisfies \eqref{equation:li-filtration-property}, and $\Hgy(\Cpx(V))_{(p)}$ contains $\Hgy(\Cpx(V)(p)) = \Hgy(\Cpx(V))(p)$ for all $p$ in $\Z$, so by \cite[Theorem 4.14]{li2004vertex}, there is an inclusion 
\[ \Hgy(\Cpx(V))_{p} \subseteq \Hgy(\Cpx(V))_{(p)}. \]

%
%Because the grading and the vertex algebra structure on the cohomology is induced by the grading and the structure on the complex, a direct computation shows that the induced filtration is the increasing Li filtration of $\Hgy(\Cpx(V))$.

The following statement gives a sufficient condition on $X$ to compute $\gr \Hgy(\Cpx(V))$. Variations of this statement appear frequently in the literature \cite{arakawa2015associated,arakawa2015localization,kuwabara2021,,arakawa2024arc,genra2025reduction, arakawa2026faithful}. We give the proof in our particular setting.

\begin{proposition} \label{proposition:reduction-quantisation}
	We make the following technical assumption:
	\begin{enumerate}
		\item the moment map $\mu_x\colon X \rightarrow (\g_{\geqslant 1})^\vee$ is a smooth scheme homomorphism,
		
		\item there is an affine subscheme $S$ of $(\mu_X)^{-1}(0)$ such that the action of $G_{\geqslant 1}$ induces a scheme isomorphism
		\[ G_{\geqslant 1} \times S \cong (\mu_X)^{-1}(0). \]
	\end{enumerate}
	Then, the following properties hold:
	\begin{enumerate}
		\item there are scheme isomorphisms
		\[ \J_\infty S \cong \J_\infty (\mu_X)^{-1}(0) \git \J_\infty G_{\geqslant 1} \cong \J_\infty ((\mu_X)^{-1}(0) \git G_{\geqslant 1}), \]
		
		\item there is the following cohomology vanishing:
		\[ \Hgy^n(\Cpx(V)) = 0 \quad \text{for} \quad n \neq 0, \]
		
		\item the natural homomorphism
		\[ \gr \Hgy^0(\Cpx(V)) \longrightarrow \Hgy^0(\gr \Cpx(V))  \]
		is a vertex Poisson algebra isomorphism. Moreover, both are explicitly given by the arc space of the Hamiltonian reduction of $X$:
		\[ \gr \Hgy^0(\Cpx(V)) \cong \Hgy^0(\gr \Cpx(V)) \cong \J_\infty \C[(\mu_X)^{-1}(0) \git G_{\geqslant 1}].  \]
	\end{enumerate}
\end{proposition}

\begin{proof}
	The first statement is given by applying the functor $\J_\infty$ to the isomorphism $G_{\geqslant 1} \times S \cong (\mu_X)^{-1}(0)$.
	
	The increasing Li filtration on $\Cpx(V)$ is bounded from below and exhaustive. By \cite[Theorem 5.5.1]{weibel1994introduction}, it gives rise to a convergent spectral sequence $\{E_n\}_{n \geqslant 0}$ whose first pages and the infinity page are:
	\[ \begin{aligned}
		E_0^{p, q} & = \gr_p \Cpx(V)^{\mathrm{ch} = p + q}, \\
		E_1^{p, q} & = \Hgy^{p + q}(\gr_p \Cpx(V)), \\
		E_\infty^{p, q} & = \gr_p \Hgy^{P + q}(\Cpx(V)),
	\end{aligned} \]
	for $p$ and $q$ in $\Z$.
	
	The $0$-page is a vertex Poisson analogue of the BRST complex:
	\[ \gr \Cpx(V) = \C[\J_\infty X] \otimes \J_\infty \operatorname{Alt}(\g_{\geqslant 1} \oplus (\g_{\geqslant 1})^\vee). \]
	See \cite[Section 2.3.3]{arakawa2015localization} for precise definitions.
	
	Because $\mu_X$ is assumed to be smooth, the induced map $\J_\infty \mu_X \colon \J_\infty X \rightarrow \J_\infty (\g_{\geqslant 1})^\vee$ is flat by \cite[Remark 3.4.2]{genra2025reduction}. Moreover, we have assumed the isomorphism $G_{\geqslant 1} \times S \cong (\mu_X)^{-1}(0)$. Therefore, we can apply \cite[Theorem 2.3.3.1]{arakawa2015localization} and we get that
	\[ \Hgy^{n}(\gr\Cpx(V)) = 0 \quad \text{for} \quad n \neq 0. \]
	
	Therefore the $1$-page $E_1^{p, q}$ is zero outside the diagonal $p = - q$. Therefore, we have isomorphisms
	\[ E_r \cong E_{r+1} \cong E_\infty \quad \text{for} \quad r \geqslant 1, \]
	and in particular
	\[ \gr \Hgy^n(\Cpx(V)) \cong  \Hgy^n(\gr \Cpx(V)) \]
	for all $n$. But these two cohomologies are zero for $n \neq 0$. This gives the vanishing of $\Hgy^n(\Cpx(V))$ (because the filtration is bounded from below and exhaustive) and the vertex Poisson isomorphism
	\[ \gr \Hgy^0(\Cpx(V)) \cong \Hgy^0(\gr \Cpx(V)). \]
\end{proof}

\section{Decomposition theorem in the chiral setting} \label{section:chiral-setting}

\subsection{Affine W-algebra for the minimal orbit}

Affine W-algebras are defined by the following chiral Hamiltonian reduction process \cite{feigin1990quantization,kac2003quantum,kac2004quantum}. See also \cite[Section 4]{genra2025reduction} for a detailed survey of the construction. Let us define a chiral comoment map
\[ \widehat\mu \colon \V(\g_{\geqslant 1}) \longrightarrow \V^k(\g) \]
by mapping a generator $a$ in $\g_{\geqslant 1}$ to $a + \chi(a) 1$. It is a graded homomorphism because: $V(\g_{\geqslant 1}) \subseteq \V^k(\g)$ is a graded inclusion, and $\g_{\geqslant 1} = \g_1 \subseteq \V^k(\g)(0)$.

\begin{theorem}[{\cite[Theorem 4.1]{kac2004quantum}}]
	The cohomology of $\Cpx(\V^k(\g))$ is concentrated in degree $0$.
\end{theorem}

The degree $0$ of the cohomology is, by definition, the \emph{affine W-algebra} associated to $f$:
\[ \W^k(\g, f) \coloneqq \Hgy^0\big(\Cpx(\V^k(\g))\big). \]
According to \Cref{proposition:brst-cohomology}, the W-algebra is a graded vertex algebra whose grading is induced by \eqref{equation:grading-universal-affine}.

\begin{remark}
	Similarly to the finite W-algebra case, the fact that the grading is supported on $\Z_{\geqslant 0}$ comes from the fact that the good grading on $\g$ is supported on $\{-1, 0, 1\}$ (see \Cref{remark:nonnegative-filtration}).
\end{remark}

\subsection{The W-algebra is a vertex subalgebra of the localisation} \label{subsection:affine-w-algebra-coset}

We state a chiralisation of \Cref{theorem:slodowy-git}, that is the reconstruction of the affine W-algebra as a subalgebra of $\V^k(\g)[e^{-1}]$. 

The Lie algebra $\J_\infty(\kfk \ltimes \g_{\geqslant 1}) = (\kfk \ltimes \g_{\geqslant 1})[\![t]\!]$, as a subalgebra of $\J_\infty \g$, acts on $\V^k(\g)$ and the localisation $\V^k(\g)[e^{-1}]$. 

The Lie algebra $\J_\infty(\kfk \ltimes \g_{\geqslant 1})$ inherits a $\Z_{\leqslant 0}$-grading from the opposite good grading $\g = \bigoplus_{i \in \Z} \g_{-i}$ and by imposing that the formal variable $t$ is of degree $0$ (see \Cref{example:grading-jet-groups}). The properties of the Lie filtration (\Cref{theorem:li-filtration}) give that for $x$ in $\g_i \subseteq \V^k(\g)(1 - i)$, and for nonnegative integers $n$ and $p$,
\[ \begin{multlined}
	x t^n \cdot \big(\V^k(\g)[e^{-1}]\big)_p = x_{(n)} \big(\V^k(\g)[e^{-1}]\big)_p \\ \subseteq \big(\V^k(\g)[e^{-1}]\big)_{1 - i + p -1} = \big(\V^k(\g)[e^{-1}]\big)_{p - i}.
\end{multlined} \]
Therefore, $\V^k(\g)[e^{-1}]$ is a filtered $\J_\infty(\kfk \ltimes \g_{\geqslant 1})$-module in the sense of \Cref{subsection:filtered-representation}.

Because the action of $x t^n$ is given by the $n$-product $x_{(n)}$, we have:
\[ \begin{aligned}
	& \big(\V^k(\g)[e^{-1}] \big)^{\J_\infty(\kfk \ltimes \g_{\geqslant 1})} \\
	= & \big\{ F \in \V^k(\g)[e^{-1}] \quad | \quad \text{for all} ~  x \in \kfk \ltimes \g_{\geqslant 1} ~\text{and} ~n \geqslant 0,  ~  x_{(n)} F = 0 \big\} \\
	= & \big\{ F \in \V^k(\g)[e^{-1}] \quad | \quad \text{for all} ~  x \in \V^k(\kfk \ltimes \g_{\geqslant 1}),  ~   [x(z_1), F(z_2)] = 0 \big\}. 
\end{aligned} \]
Because it is defined as subalgebra commuting with a graded vertex subalgebra, $\big(\V^k(\g)[e^{-1}] \big)^{\J_\infty(\kfk \ltimes \g_{\geqslant 1})}$ is a graded vertex subalgebra of $\V^k(\g)[e^{-1}]$.

\begin{theorem} \label{theorem:affine-w-algebra-coset}
	There is a graded vertex algebra isomorphism
	\[ \big(\V^k(\g)[e^{-1}]\big)^{\J_\infty(\kfk \ltimes \g_{\geqslant 1})} \cong \W^k(\g, f) \]
\end{theorem}

To prove this theorem, we start by computing the associated graded algebra of $\big(\V^k(\g)[e^{-1}]\big)^{\J_\infty(\kfk \ltimes \g_{\geqslant 1})}$ for the filtration induced by the increasing Li filtration on $\V^k(\g)[e^{-1}]$.

\begin{proposition} \label{proposition:gr-invariant-borisov-localisation}
	The natural vertex Poisson algebra homomorphism
	\[ \gr\big((\V^k(\g)[e^{-1}])^{\J_\infty(\kfk \ltimes \g_{\geqslant 1})}\big) \longrightarrow \big(\gr(\V^k(\g)[e^{-1}])\big)^{\J_\infty(\kfk \ltimes \g_{\geqslant 1})} \]
	is an isomorphism. Therefore, we have a vertex Poisson algebra isomorphism:
	\[ \gr\big((\V^k(\g)[e^{-1}])^{\J_\infty(\kfk \ltimes \g_{\geqslant 1})}\big) \cong \J_\infty\C[\g^\vee[e^{-1}] \git (K \ltimes G_{\geqslant 1})].  \]
\end{proposition}

\begin{proof}
	For an algebraic group $N$, and for any $m$ in $\Z_{\geqslant 0}$, there is a semidirect product decomposition
	\[ \J_m N = N \ltimes (\J_m N)_1, \]
	where $(\J_m N)_1$ is the normal subgroup given by the kernel of the natural projection $\J_m N \rightarrow N$. The group $(\J_m N)_1$ is unipotent for finite $m$ and its Lie algebra is $t \n[t] / (t^{m+1})$.
	
	Therefore, the arc group $\J_\infty K = \J_\infty \Gm \times \J_\infty \Ga{}^{d-2}$ decomposes as the product of a torus and the projective limit of a system of unipotent groups:
	\[ \J_\infty K = \Gm \ltimes ((\J_\infty \Gm)_1 \times \J_\infty \Ga{}^{d-2}). \]
	
	To prove the proposition, we can first deal with the action of the normal subgroups: $((\J_\infty \Gm)_1 \times \J_\infty \Ga{}^{d-2}) \ltimes \J_\infty G_{\geqslant 1}$ in $\J_\infty K \ltimes \J_\infty G_{\geqslant 1}$, that is the projective limit of graded unipotent groups $((\J_m \Gm)_1 \times \J_m \Ga{}^{d-2}) \ltimes \J_m G_{\geqslant 1}$. We aim to show the isomorphism
	\begin{equation} \label{equation:gr-invariant-unipotent-jet} \gr \Big(\big(\V^k(\g)[e^{-1}]\big)^{(t \g_{\mathrm m}[\![t]\!] \oplus \g_{\mathrm a}{}^{d-2} [\![t]\!]) \ltimes \g_{\geqslant 1}[\![t]\!]}\Big) \cong \big(\gr \big(\V^k(\g)[e^{-1}]\big)\big)^{(t \g_{\mathrm m}[\![t]\!] \oplus \g_{\mathrm a}{}^{d-2} [\![t]\!]) \ltimes \g_{\geqslant 1}[\![t]\!]}. \end{equation}
	To do so, we use \Cref{proposition:gr-invariant-pro}. If we apply the arc functor $\J_\infty$ to the isomorphism of \Cref{proposition:dilation-gan-ginzburg}, we have the $((\J_\infty \Gm)_1 \times \J_\infty \Ga{}^{d-2}) \ltimes \J_\infty G_{\geqslant 1}$-equivariant isomorphism
	\begin{multline*} \C[((\J_\infty \Gm)_1 \times \J_\infty \Ga{}^{d-2}) \ltimes \J_\infty G_{\geqslant 1}] \otimes \C[\Gm \times \J_\infty S_f] \\ \cong \J_\infty (\Sym \g)[e^{-1}] \cong \gr (\V^k(\g)[e^{-1}]).  \end{multline*}
	The assumption of \Cref{proposition:gr-invariant-pro} about the Lie algebra gradings has been checked at the beginning of \Cref{subsection:affine-w-algebra-coset}. Therefore, we can apply the proposition and we get the isomorphism \eqref{equation:gr-invariant-unipotent-jet}.
	
	The isomorphism \eqref{equation:gr-invariant-unipotent-jet} is $\Gm$-equivariant, so we get the isomorphism
	\[ \Big(\big(\gr \big(\V^k(\g)[e^{-1}]\big)\big)^{(t \g_{\mathrm m}[t] \oplus \g_{\mathrm a}{}^{d-2} [t]) \ltimes \g_{\geqslant 1}[t]}\Big)^{\Gm} \cong \big(\gr \big(\V^k(\g)[e^{-1}]\big)\big)^{(\kfk \ltimes \g_{\geqslant 1})[t]} \]
	between the $\Gm$-invariant subalgebras. By applying \Cref{proposition:weight-filtration} on the left-hand side, we get the desired isomorphism:
	\[ \gr\big((\V^k(\g)[e^{-1}])^{\J_\infty(\kfk \ltimes \g_{\geqslant 1})}\big) \cong \big(\gr(\V^k(\g)[e^{-1}])\big)^{\J_\infty(\kfk \ltimes \g_{\geqslant 1})}. \]
	The first part of the statement is proved.
	
	We immediately deduce the isomorphisms 
	\[ \begin{aligned}
		\gr\big((\V^k(\g)[e^{-1}])^{\J_\infty(\kfk \ltimes \g_{\geqslant 1})}\big) & \cong \big(\J_\infty(\Sym \g)[e^{-1}]\big)^{\J_\infty(\kfk \ltimes \g_{\geqslant 1})}  \\
		& \cong  \big(\J_\infty(\Sym \g)[e^{-1}]\big)^{\J_\infty(K \ltimes G_{\geqslant 1})} \\
		& \cong \C[\g^\vee[e^{-1}] \git \J_\infty(K \ltimes G_{\geqslant 1})]
	\end{aligned} \]	
	by applying \Cref{lemma:coholology-pro-algebra}. As a consequence of \Cref{proposition:dilation-gan-ginzburg}, we have scheme isomorphisms
	\[ \J_\infty S_f \cong \g^\vee[e^{-1}] \git \J_\infty(K \ltimes G_{\geqslant 1}) \cong \J_\infty\big(\g^\vee[e^{-1}] \git (K \ltimes G_{\geqslant 1})\big).  \]
	This gives the second part of the statement.
\end{proof}

The following lemma gives a natural vertex algebra homomorphism between the invariant subalgebra $\big(\V^k(\g)[e^{-1}]\big)^{\J_\infty(\kfk \ltimes \g_{\geqslant 1})}$ and the chiral Hamiltonian reduction of the localisation $\V^k(\g)[e^{-1}]$. 

\begin{lemma} \label{lemma:coset-to-cohomology}
	The image of the graded vertex algebra homomorphism
	\[ \begin{aligned}
		(\V^k(\g)[e^{-1}])^{\J_\infty(\kfk \ltimes \g_{\geqslant 1})} & \longrightarrow \Cpx(\V^k(\g)[e^{-1}] = \V^k(\g)[e^{-1}] \otimes \F(\g_{\geqslant 1} \oplus (\g_{\geqslant 1})^\vee), \\
		v & \longmapsto v \otimes 1 
	\end{aligned} \]
	is in the kernel of the differential $\widehat\mu_{e^{-1}}(Q)_{(0)}$.
\end{lemma}

\begin{proof}
	Let $v$ be in $\V^k(\g)[e^{-1}]$. By using the explicit formula for $Q$ and Borcherds identity \eqref{equation:borcherds}, we get:
	\[ \widehat\mu_{e^{-1}}(Q)_{(0)} (v \otimes 1) = \sum_{i=1}^{d-1} (a_i \otimes \phi(a_i^\vee))_{(0)} (v \otimes 1) \]
	because $v$ commutes with the generators of $\F(\g_{\geqslant 1} \oplus (\g_{\geqslant 1})^\vee)$. Then we use the Borcherds identity \eqref{equation:borcherds}:
	\[ \begin{multlined}
		(a_i \otimes \phi(a_i^\vee))_{(0)} v = \sum_{n \geqslant 0} (a_i \otimes 1)_{(-n - 1)} \big((1 \otimes \phi(a_i^\vee))_{(n)} (v \otimes 1)\big) \\
		+ \sum_{n \geqslant 0} (1 \otimes \phi(a_i^\vee))_{(-n - 1)} \big((a_i{}_{(n)} v) \otimes 1\big).
	\end{multlined}  \] 
	In the previous formula, the first sum vanishes because $v$ commutes with the generators of $\F(\g_{\geqslant 1} \oplus (\g_{\geqslant 1})^\vee)$. By using the equality $a_i{}_{(n)} v = a_i t^n \cdot v$, we finally get:
	\[ \widehat\mu_{e^{-1}}(Q)_{(0)} (v \otimes 1) =  \sum_{n \geqslant 0} (1 \otimes \phi(a_i^\vee))_{(-n - 1)} \big((a_i t^n \cdot v) \otimes 1\big). \]
	By definition of being invariant, this quantity is zero if $v$ belongs to the invariant subalgebra $(\V^k(\g)[e^{-1}])^{\J_\infty(\kfk \ltimes \g_{\geqslant 1})}$.
\end{proof}

Now we can prove the main theorem.

\begin{proof}[Proof of \Cref{theorem:affine-w-algebra-coset}]
	Let us show that the localisation at $e$ does not modify the cohomology. By composing $\widehat{\mu}$ with the localisation homomorphism, we get another chiral comoment map
	\[ \widehat\mu_{e^{-1}} \colon \V(\g_{\geqslant 1}) \longrightarrow \V^k(\g)[e^{-1}]. \]
	Because the BRST complex construction is functorial, there is a natural cochain complex and vertex homomorphism
	\[ \Cpx(\V^k(\g)) \longrightarrow \Cpx(\V^k(\g)[e^{-1}]) = \Cpx(\V^k(\g))[e^{-1}] \]
	that is also a localisation homomorphism because tensoring with the fermionic vertex superalgebra commutes with localisation. Then we get a vertex superalgebra homomorphism by taking cohomologies:
	\[ \W^k(\g, f) = \Hgy\big(\Cpx(\V^k(\g))\big) \longrightarrow  \Hgy\big(\Cpx(\V^k(\g)[e^{-1}])\big). \]
	Let us prove that it is a vertex algebra isomorphism by using \Cref{proposition:reduction-quantisation}.

	It is straightforward to check that the moment maps $\mu\colon \g^\vee \rightarrow (\g_{\geqslant 1})^\vee$ and $\mu_{e^{-1}}\colon \g^\vee[e^{-1}] \rightarrow (\g_{\geqslant 1})^\vee$ are smooth by checking that their differential are surjective at any point. Moreover, $\J_\infty \mu^*$ coincides with $\gr \widehat{\mu}$ up to the isomorphism $\J_\infty(\Sym \g) \cong \gr \V^k(\g)$. \Cref{theorem:gan-ginzburg-isomorphism} and \Cref{proposition:hamiltonian-reduction-open} give the other assumptions that are necessary to apply \Cref{proposition:reduction-quantisation}. Therefore, both cohomologies of $\Cpx(\V^k(\g))$ and $\Cpx(\V^k(\g))[e^{-1}]$ are concentrated in degree $0$ and we have the following commutative diagram of vertex Poisson algebra homomorphisms:
	\[ \begin{tikzcd}
		\gr \Hgy^0(\Cpx(\V^k(g))) \arrow[r] \arrow[d, "\sim"'{sloped}] & \gr \Hgy^0\big(\Cpx(\V^k(\g))[e^{-1}]\big) \arrow[d, "\sim"{sloped}] \\
		\J_\infty \C[\mu^{-1}(0) \git G_{\geqslant 1}] \arrow[r, "\sim"] & \J_\infty \C[(\mu_{e^{-1}})^{-1}(0) \git G_{\geqslant 1}],
	\end{tikzcd} \]
	where the vertical arrows are isomorphisms. In fact, the lower horizontal arrow is an isomorphism by applying the arc functor $\J_\infty$ to the isomorphism from \Cref{proposition:hamiltonian-reduction-open}. Therefore, we get a vertex Poisson algebra isomorphism
	\[ \gr \W^k(\g, f) \cong \gr \Hgy^0\big(\Cpx(\V^k(\g)[e^{-1}]\big) \]
	that is given by the upper arrow. Since the filtrations on the cohomology vertex algebras are induced by the Li filtrations on the complexes, these filtrations are bounded from below and exhaustive. By \Cref{lemma:quantisation-iso}, we get the vertex algebra isomorphism:
	\[ \W^k(\g, f) \cong \Hgy^0\big(\Cpx(\V^k(\g)[e^{-1}]\big) = \Hgy\big(\Cpx(\V^k(\g)[e^{-1}]\big). \]
	
	As a consequence of \Cref{lemma:coset-to-cohomology}, we have a vertex algebra homomorphism
	\[ \widehat\lambda \colon (\V^k(\g)[e^{-1}])^{\J_\infty(\kfk \ltimes \g_{\geqslant 1})} \longrightarrow \Hgy^0\big(\Cpx(\V^k(\g)[e^{-1}]\big) \cong \W^k(\g, f). \]
	Both sides are equipped with exhaustive filtrations that are bounded from below, so it is an isomorphism if $\gr \widehat\lambda$ is an isomorphism, according to \Cref{lemma:quantisation-iso}. 
	
	By \Cref{proposition:gr-invariant-borisov-localisation} and \Cref{proposition:reduction-quantisation} (see the first part of this proof), we have the algebra isomorphisms
	\[ \begin{multlined}
		\gr \big(\big(\V^k(\g)[e^{-1}]\big)^{\J_\infty(\kfk \ltimes \g_{\geqslant 1})}\big) \cong \J_\infty \C[\g^\vee[e^{-1}] \git (K \ltimes G_{\geqslant 1})] \\ \text{and} \quad \gr \W^k(\g, f) \cong \J_\infty \C[\mu^{-1}(0) \git G_{\geqslant 1}]. 
	\end{multlined}\]
	Under these isomorphisms, $\gr \widehat\lambda$ is conjugate to the isomorphism
	\[ \J_\infty \C[\g^\vee[e^{-1}] \git (K \ltimes G_{\geqslant 1})] \cong  \J_\infty \C[\mu^{-1}(0) \git G_{\geqslant 1}] \]
	given by \Cref{theorem:slodowy-git}. Therefore $\gr \widehat\lambda$, and then $\widehat\lambda$, are isomorphisms.
\end{proof}

\subsection{Decomposition theorem for affine W-algebras}

\begin{lemma}
	The two vertex subalgebras 
	\[ \V^k(\kfk \ltimes \g_{\geqslant 1})[e^{-1}] \quad \text{and} \quad \big(\V^k(\g)[e^{-1}]\big)^{\J_\infty(\kfk \ltimes \g_{\geqslant 1})} \cong \W^k(\g, f) \]
	of $\V^k(\g)[e^{-1}]$ commute with each other.
\end{lemma}

\begin{proof}
	The lemma follows immediately from the following equalities:
	\[ \begin{aligned}
		& \big(\V^k(\g)[e^{-1}] \big)^{\J_\infty(\kfk \ltimes \g_{\geqslant 1})} \\
		= & \big\{ F \in \V^k(\g)[e^{-1}] \quad | \quad \text{for all} ~  x \in \kfk \ltimes \g_{\geqslant 1} ~\text{and} ~n \geqslant 0,  ~   xt^n \cdot F = x_{(n)} F = 0 \big\} \\
		= & \big\{ F \in \V^k(\g)[e^{-1}] \quad | \quad \text{for all} ~  x \in \V^k(\kfk \ltimes \g_{\geqslant 1})[e^{-1}],  ~   [x(z_1), F(z_2)] = 0 \big\}. 
	\end{aligned} \]
	The second equality follows because commuting with $\V^k(\kfk \ltimes \g_{\geqslant 1})[e^{-1}]$ is equivalent to commuting with $\kfk \ltimes \g_{\geqslant 1}$, since this space generates the non-localised algebra, and commuting with $e(z_1)$ implies commuting with $e^{-1}(z_1) = (e(z_1))^{-1}$ by \cite[Formula (3.1)]{arakawa2026faithful}.
	
\end{proof}

Consequently, there is a natural vertex algebra homomorphism
\[ \widehat{\delta} \colon \V^k(\kfk \ltimes \g_{\geqslant 1})[e^{-1}] \otimes \big(\V^k(\g)[e^{-1}]\big)^{\J_\infty(\kfk \ltimes \g_{\geqslant 1})} \longrightarrow \V^k(\g)[e^{-1}] \]
given by the $(-1)$-product. This provides the decomposition isomorphism for the affine W-algebra.

\begin{theorem} \label{theorem:decomposition-chiral}
	The vertex algebra homomorphism
	\[ \widehat{\delta} \colon \V^k(\kfk \ltimes \g_{\geqslant 1})[e^{-1}] \otimes \big(\V^k(\g)[e^{-1}]\big)^{\J_\infty(\kfk \ltimes \g_{\geqslant 1})} \longrightarrow \V^k(\g)[e^{-1}] \]
	is an isomorphism.
\end{theorem}

\begin{proof}	
	Both sides are equipped with increasing Li filtrations. Because these filtrations on both sides are bounded from below and exhaustive, we just need to prove that
	\[ \gr \widehat\delta \colon \gr \big(\V^k(\kfk \ltimes \g_{\geqslant 1})[e^{-1}]\big) \otimes \gr\big(\big(\V^k(\g)[e^{-1}]\big)^{\J_\infty(\kfk \ltimes \g_{\geqslant 1})}\big) \longrightarrow \gr \big(\V^k(\g)[e^{-1}]\big) \]
	is an isomorphism to apply \Cref{lemma:quantisation-iso}.
	
	By applying \Cref{proposition:gr-borisov-localisation}, we obtain the vertex Poisson algebra isomorphisms:
	\[ \gr\big(\V^k(\kfk \ltimes \g_{\geqslant 1})[e^{-1}]\big) \cong \J_\infty\Sym(\kfk \ltimes \g_{\geqslant 1})[e^{-1}]  \]
	and 
	\[ \gr\big(\V^k(\g)[e^{-1}]\big) \cong \J_\infty \Sym(\g)[e^{-1}]. \]
	
	Moreover, \Cref{proposition:gr-invariant-borisov-localisation} gives the vertex Poisson algebra isomorphism
	\[ \gr\big((\V^k(\g)[e^{-1}])^{\J_\infty(\kfk \ltimes \g_{\geqslant 1})}\big) \cong \J_\infty\C[\g^\vee[e^{-1}] \git(K \ltimes G_{\geqslant 1})].  \]
	
	Under these isomorphisms, $\gr \widehat{\delta}$ is conjugate to the vertex Poisson algebra isomorphism 
	\[ \J_\infty\delta^* \colon \J_\infty \C[(\kfk \ltimes \g_{\geqslant 1})^\vee][e^{-1}] \otimes \J_\infty\C[\g^\vee[e^{-1}] \git (K \ltimes G_{\geqslant 1})]  \overset{\sim}{\longrightarrow} \J_\infty \C[\g^\vee[e^{-1}]]  \]
	provided by \Cref{theorem:decomposition-slodowy}. This concludes the proof.
\end{proof}

\subsection{Inverse Hamiltonian reduction for the minimal W-algebra}

If we compose the localisation embedding $\V^k(\g) \hookrightarrow \V^k(\g)[e^{-1}]$ with the inverse $\widehat{\delta}^{-1}$ of the decomposition isomorphism (\Cref{theorem:decomposition-chiral}), we get an embedding of graded vertex algebras:
\[ \V^k(\g) \hooklongrightarrow \V^k(\kfk \ltimes \g_{\geqslant 1}) \otimes \W^k(\g, f). \]
This is similar to the inverse reduction embedding \eqref{equation:inverse-reduction}. To check whether it is an inverse reduction, we need to determine if the factor $\V^k(\kfk \ltimes \g_{\geqslant 1})$ is a \emph{localised Weyl vertex algebra}. 

We recall that the \emph{Weyl vertex algebra} of rank $r$ (also called \emph{$\beta\gamma$-system} in the physics literature) is a vertex algebra freely generated by $2 r$ generators $\beta_i$ and $\gamma_i$, where $1 \leqslant i \leqslant r$, whose only non-trivial OPE are given by:
\[ \beta_i(z_1) \gamma_i(z_2) \sim \frac{1}{z_1 - z_2}, \quad \text{for} \quad 1 \leqslant i \leqslant r. \]
We denote this vertex algebra by $\mathrm{Weyl}_r$. It can be equipped with a vertex algebra grading by imposing that the $\beta_i$ are in degree $1$ and the $\gamma_i$ are in degree $0$. In particular, the localisation $\mathrm{Weyl}_r[(\gamma_1)^{-1}]$ is well-defined. See for example \cite[Section 1.8]{desole2006finite}, where it corresponds to a purely even fermionic vertex algebra associated with the canonical symplectic form on $\C^r \times \C^r$. 

For $2 \leqslant i \leqslant d - 1$, we introduce the notations
\[ a_i \coloneqq x_{i + 1, 2} \quad \text{and} \quad b_i \coloneqq x_{1, i + 1}. \]
Then $\kfk \ltimes \g_{\geqslant 1}$ is spanned by $H$, $e$, $a_i$ and $b_i$ for $2 \leqslant i \leqslant d - 1$, with the only non-trivial bracket relations given by:
\[ [H, e] = e, \quad [H, b_i] = b_i, \quad [b_i, a_i] =e. \]

\begin{proposition} \label{proposition:relations-dual-w-algebra}
	The vertex algebra $\V^k(\kfk \ltimes \g_{\geqslant 1})[e^{-1}]$ is strongly generated by $H$, $e$, $e^{-1}$, $a_i$ and $b_i$ for $2 \leqslant i \leqslant d - 1$, with the only non-trivial OPE given by:
	\[ \begin{aligned}
		H(z_1) H(z_2) & \sim \frac{1}{(z_1 - z_2)^2} \frac{d-1}{d} k \id, \\
		H(z_1) e(z_2) & \sim \frac{1}{z_1 - z_2} e(z_2), \\
		H(z_1) e^{-1}(z_2) & \sim \frac{1}{z_1 - z_2} \big(- e(z_2)\big), \\
		H(z_1) b_i(z_2) & \sim \frac{1}{(z_1 - z_2)^1} b_i(z_2), \\
		b_i(z_1) a_i(z_2) & \sim \frac{1}{z_1 - z_2} e(z_2).
	\end{aligned}\]
\end{proposition}

\begin{proof}
	The only relation that must be computed is between $H$ and $e^{-1}$. This can be done with computations analogous to the proof of \Cref{proposition:relations-localisation-sl2}.
\end{proof}

From this, we deduce that this vertex algebra is a localised Weyl vertex algebra.

\begin{theorem} \label{theorem:commutant-localisation-weyl}
	The assignment 
	\[ \begin{aligned}
		\mathrm{Weyl}_{d-1} & \longrightarrow  \V^k(\kfk \ltimes \g_{\geqslant 1})[e^{-1}], \\
		\beta_1 & \longmapsto  \NO{e^{-1} H} + (3 - \frac{d - 1}{4 d} k) e^{-1} e^{-1} \partial e, \\
		\gamma_1 & \longmapsto e, \\
		\beta_i & \longmapsto - a_i, \\
		\gamma_i & \longmapsto e^{-1} b_i \quad \text{for} \quad 2 \leqslant i \leqslant d - 1,
	\end{aligned}  \]
	is a well-defined graded vertex algebra homomorphism which induces an isomorphism
	\[ \mathrm{Weyl}_{d-1}[(\gamma_1)^{-1}] \cong \V^k(\kfk \ltimes \g_{\geqslant 1})[e^{-1}]. \]
\end{theorem}

\begin{proof}
	The mapping is clearly degree preserving because $H$, the $a_i$ and $\partial e$ are of degree $1$, whereas $e$, $e^{-1}$ and the $b_i$ are of degree $0$.
	
	By using the Borcherds identity \eqref{equation:borcherds}, we get:
	\begin{equation} \label{equation:ope-quasi-beta-gamma} \begin{aligned}
		\NO{e^{-1} H}(z_1) e(z_2) & \sim \frac{1}{z_1 - z_2} \id, \\ 
		\NO{e^{-1} H}(z_1) (e^{-1} b_i)(z_2) & \sim 0, \\
		(e^{-1} b_i) (z_1) a_i(z_2) & \sim \frac{1}{z_1 - z_2} \id.
	\end{aligned} \end{equation}
	
	However, the OPE of $\NO{e^{-1}H}$ with itself is non-trivial:
	\[ \begin{multlined}
		\NO{e^{-1}H}(z_1) \NO{e^{-1}H}(z_2) \sim \frac{1}{z_1 - z_2} \big(12 - \frac{d-1}{d} k \big) (e^{-1} e^{-1} e^{-1} \partial e)(z_2) \\
		 + \frac{1}{(z_1 - z_2)^2} \big(\frac{d-1}{d} k - 12\big) e^{-1} e^{-1}.
	\end{multlined} \]
	By replacing $\NO{e^{-1}H}$ by the corrected term 
	\[ c \coloneqq  \NO{e^{-1}H} + (3 - \frac{d - 1}{4 d} k) e^{-1} e^{-1} \partial e, \]
	we get the same OPE as \eqref{equation:ope-quasi-beta-gamma} because $e$ and $e^{-1}$ commute with the $b_i$ and $a_i$. In addition, we get a trivial relation 
	\[ c(z_1) c(z_2) \sim 0. \]
	
	We have checked that the elements $c$, $e$, $e^{-1} b_i$ and $a_i$ have the OPE of the Weyl vertex algebra generators $\gamma_j$ and $\beta_j$. Therefore, we have a well-defined vertex algebra homomorphism.  
	
	To prove that the induced vertex algebra homomorphism \[ \mathrm{Weyl}_{d-1}[(\gamma_1)^{-1}] \longrightarrow \V^k(\kfk \ltimes \g_{\geqslant 1})[e^{-1}], \]
	is invertible, we can use the increasing Li filtrations on both sides. The associated graded algebra of $\mathrm{Weyl}_{d-1}$ is obtained by applying the arc functor $\J_\infty$ to a polynomial algebra with $2 (d - 1)$ generators denoted by $\beta_i$ and $\gamma_i$ (the algebra of polynomial functions over a symplectic vector space $\C^{d - 1} \times \C^{d - 1}$):
	\[ \gr \mathrm{Weyl}_{d-1} \cong \J_\infty \C[\gamma_i, \beta_i]_{1 \leqslant i \leqslant d - 1}. \]
	After taking $\gr$, the induced homomorphism is given by
	\[ \begin{aligned}
		\gr \mathrm{Weyl}_{d-1}[(\gamma_1)^{-1}] & \longrightarrow  \gr \V^k(\kfk \ltimes \g_{\geqslant 1})[e^{-1}], \\
		\beta_1 & \longmapsto  e^{-1} H, \\
		\gamma_1 & \longmapsto e, \\
		\beta_i & \longmapsto - a_i, \\
		\gamma_i & \longmapsto e^{-1} b_i \quad \text{for} \quad 2 \leqslant i \leqslant d - 1.
	\end{aligned}  \] 
	It is clearly an isomorphism whose inverse map is given by:
	\[  H \longmapsto \beta_1 \gamma_1, \quad e \longmapsto \gamma_1, \quad a_i \longmapsto - \beta_i, \quad b_i \longmapsto \gamma_1 \gamma_i.  \]
	Therefore, \Cref{lemma:quantisation-iso} implies that the vertex algebra homomorphism we have defined is an isomorphism.
\end{proof}

Therefore, we retrieve the inverse reduction for the minimal W-algebra in type~A as a corollary of \Cref{theorem:decomposition-chiral}.

\begin{theorem}[\cite{semikhatov1994inverting,adamovic2019realizations,adamovic2021realisation,fehily2024inverse}]
	There is a graded vertex algebra embedding
	\[ \V^k(\g) \hooklongrightarrow \mathrm{Weyl}_{d-1}[(\gamma_1)^{-1}] \otimes \W^k(\g, f). \]
\end{theorem}

\section{The $\sl_2$-example in the chiral setting, at any level} \label{section:sl2-example}

\subsection{Definition of the localisation}

The universal affine vertex algebra $\V^k(\sl_2)$ is strongly generated by three vectors $e$, $h$ and $f$ whose only non-trivial OPE are given by
\[ \begin{multlined}
	h(z_1) h(z_2) \sim \frac{1}{(z_1 - z_2)^2} 2 k \id, \quad h(z_1) e(z_2) \sim \frac{1}{z_1 - z_2} 2 e(z_2) , \\ e(z_1) f(z_2) \sim \frac{1}{z_1 - z_2} h(z_2) + \frac{1}{(z_1 - z_2)^2} k \id \\ \text{and} \quad h(z_1) f(z_2) \sim - \frac{1}{z_1 - z_2} 2 f(z_2).
\end{multlined} \]
It is equipped with the grading \eqref{equation:grading-universal-affine} such that $e$ is in degree $0$, $h$ is in degree $1$ and $f$ in degree $2$.

By the general theory of \cite{borisov2001vertex,arakawa2026faithful}, the localised graded vertex algebra $\V^k(\sl_2)[e^{-1}]$ contains an additional generator $e^{-1}$, of degree $0$. Its associated field
\[ e^{-1}(z) \coloneqq \sum_{n \in \Z} e^{-1}_{(n)} \, z^{- n -1} \]
satisfies the following properties:
\begin{enumerate}
	\item All the coefficients $e^{-1}_{(n)}$ and $e_{(n)}$, for all $n$ in $\Z$, commute with each other, that is to say
	\[ e^{-1}(z_1) e(z_2) \sim 0.  \]
	
	\item The product
	\[ e(z) e^{-1}(z) \coloneqq \sum_{n \in \Z} \Big(\sum_{i + j = n} e_{(i)} \circ e^{-1}_{(j)} \Big) z^{- n -2} \]
	is well defined, meaning that the sum $\sum_{i + j = n} e_{(i)} \circ e^{-1}_{(j)} v$ is finite for any $v$ in $\V^k(\sl_2)[e^{-1}]$. This follows because of the vertex algebra grading: $e_{(i)} v$ or $e^{-1}_{(j)} v$ is zero for big values of $i$ and $j$.
	
	\item The identity 
	\[ e^{-1}(z) e(z) = e(z) e^{-1}(z) = \id \]
	is satisfied.
\end{enumerate}

As explained in \cite[Lemma 3.1 and before Lemma 3.3]{arakawa2026faithful}, the product $\alpha \coloneq e_{(-1)}$ is an invertible endomorphism of $\V^k(\sl_2)[e^{-1}]$, but in general
\[ \alpha^{-1} = (e_{(-1)})^{-1} \neq e^{-1}_{(-1)} \coloneqq (e^{-1})_{(-1)}, \] 
see \Cref{lemma:inverse-no}.

Everything in $\V^k(\sl_2)[e^{-1}]$ can be computed by considering sums of elements of the form
\[ \alpha^{-p}(v) \quad \text{for} \quad p \in \Z_{\geqslant 0} \quad \text{and} \quad v \in \V^k(\sl_2). \]
To compute the $e^{-1}_{(n)} v$, for $n$ in $\Z$ and $v$ in $\V^k(\sl_2)$, we must use the identity
\[ e(z) e^{-1}(z) v = v. \]
If $v$ is homogeneous of degree $\Delta$, then many terms equal zero and one gets
\[ \sum_{n \leqslant \Delta - 2} \Big(\sum_{\substack{i, j \leqslant \Delta - 1 \\ i + j = n}} e_{(i)} \circ e^{-1}_{(j)} v \Big) z^{- n -2} = v. \]

If we apply this identity to the vacuum $1$, we get:
\[ \sum_{n \leqslant - 2} \Big(\sum_{\substack{i, j \leqslant - 1 \\ i + j = n}} e_{(i)} \circ e^{-1}_{(j)} 1 \Big) z^{- n -2} = 1. \]
For $n = 2$, the identification gives
\[ e_{(-1)} \circ e^{-1}_{(j)} 1 = \NO{ e e^{-1}} = 1. \]
Therefore, the equality $e^{-1}_{(-1)} = ( e_{(-1)})^{-1}$ holds on the commutative vertex subalgebra $\big(\V^k(\sl_2)[e^{-1}]\big)_0$ which is generated by $e$ and $e^{-1}$ (it is the lowest part of the increasing Li filtration): $e^{-1}$ is the inverse of $e$ in the usual associative sense.

The translation operator $\partial$ satisfies a Leibniz rule for the normally ordered product:
\[ \partial(\NO{ab}) = \NO{\partial(a) b} + \NO{a \partial(b)}. \]
Therefore, we get
\[ \partial(e^{-1}) = - \NO{e \NO{e \partial(e)}} = - e^{-2} \partial e, \]
where we can drop the normally ordered product because everything commutes.

\emph{In general, if an expression only involves $e$, $e^{-1}$ and their derivatives, we omit the $\NO{}$ to denote the product because everything will commute.}

\subsection{Operator product expansions} \label{subsection:ope-localised-sl2}

We compute the OPE between $e^{-1}$ and the other generators manually, by using the machinery described above.

\begin{proposition} \label{proposition:relations-localisation-sl2}
	The localised universal affine vertex algebra $\V^k(\sl_2)[e^{-1}]$ is strongly generated by four vectors $e$, $e^{-1}$, $h$ and $f$ whose only non-trivial OPE are given by
	\[ \begin{multlined}
		h(z_1) h(z_2) \sim \frac{1}{(z_1 - z_2)^2} 2 k \id, \quad h(z_1) e(z_2) \sim \frac{1}{z_1 - z_2} 2 e(z_2) , \\ e(z_1) f(z_2) \sim \frac{1}{z_1 - z_2} h(z_2) + \frac{1}{(z_1 - z_2)^2} k \id, \\
		h(z_1) f(z_2) \sim - \frac{1}{z_1 - z_2} 2 f(z_2),\quad
		e^{-1}(z_1) h(z_2) \sim  \frac{1}{z_1 - z_2} 2 e^{-1}(z_2), \\
		\text{and} \quad e^{-1}(z_1) f(z_2) \sim  \frac{1}{z_1 - z_2} \Big( 2 (k - 1) \big((e^{-1})^3 \partial e\big)(z_2) - \NO{e^{-1} e^{-1} h}(z_2)\Big) \\
		+ \frac{1}{(z_1 - z_2)^2} \Big( - (k + 2) (e^{-1})^2(z_2)\Big).
	\end{multlined} \]
	and satisfying the additional relation
	\[ e^{-1}(z) e(z) = e(z) e^{-1}(z) = \id. \]
\end{proposition}

\begin{proof}[Computation of the OPE between $e^{-1}$ and $h$]
	The operator product expansion of $e$ and $h$ is
	\[ e(z_1) h(z_2) \sim \frac{1}{z_1 - z_2} (- 2 e(z_2)), \]
	therefore $e_{(0)} h = - 2 e$ and all other products are zero. The element $h$ is of degree $1$, so the product $e^{-1}_{(n)} h$ is of degree $0 + 1 - n - 1 = -n$, that is a negative number for $n > 0$. So only the product $e^{-1}_{(0)} h$ may be non-zero. Consider the identity
	\[ \sum_{n \leqslant- 1} \Big(\sum_{\substack{i, j \leqslant 0 \\ i + j = n}} e_{(i)} \circ e^{-1}_{(j)} h \Big) z^{- n -2} = h. \]
	For $n = -1$, we get the equality
	\[ e_{(-1)} \circ e^{-1}_{(0)} h + e_{(0)} \circ e^{-1}_{(-1)} h = e_{(-1)} \circ e^{-1}_{(0)} h +  e^{-1}_{(-1)} \circ e_{(0)} h  = 0. \]
	By using the equality
	\[ e^{-1}_{(-1)} h =  e^{-1}_{(-1)} \circ e_{(0)} h = - 2 e^{-1}_{(-1)} e = - 2 , \]
	we get
	\[ e^{-1}_{(0)} h = 2 e_{(-1)} 1 = 2 e^{-1}. \]
\end{proof}

The OPE between $e^{-1}$ and $f$ is slightly more involved. We need the following lemma.

\begin{lemma} \label{lemma:inverse-no}
	The following equality holds:
	\[ (e_{(-1)})^{-1} (\NO{e^{-1} h}) - \NO{e^{-1} e^{-1} h} = 4 (e^{-1})^3 \partial e. \]
	This proves in particular that $(e_{(-1)})^{-1} \neq e^{-1}_{(-1)} = \NO{e^{-1} \bullet}$.
\end{lemma}

\begin{proof}
	Localisation imposes the equality:
	\[ e_{(0)} e^{-1}_{(-2)} (\NO{e h^{-1}}) + e_{(-1)} e^{-1}_{(-1)} (\NO{e h^{-1}}) + e_{(-2)} e^{-1}_{(0)} (\NO{e h^{-1}}) = \NO{e h^{-1}}. \]
	By the Borcherds identity \eqref{equation:borcherds}:
	\[ \begin{aligned}
		e_{(0)}(\NO{e^{-1} h^{-1}}) & = \NO{e^{-1} (e_{(0)} h)} = - 2 \NO{e^{-1}e } = - 2, \\
		e^{-1}_{(0)}(\NO{e^{-1} h^{-1}}) & = \NO{e^{-1} (e^{-1}_{(0)} h)} = 2 (e^{-1})^2.
	\end{aligned} \]
	Moreover,
	\[ \begin{aligned}
		e_{(0)} e^{-1}_{(-2)} (\NO{e h^{-1}}) & = e^{-1}_{(-2)} e_{(0)} (\NO{e h^{-1}}) = - 2 e^{-1}_{(-2)}1 = - 2 \partial(e^{-1}) = 2 (e^{-1})^2 \partial e, \\
		e_{(-2)} e^{-1}_{(0)} (\NO{e h^{-1}}) & =e_{(-2)} e^{-1}_{(0)}   (\NO{e h^{-1}}) = 2 (e^{-1})^2 \partial e.
	\end{aligned} \]
	From this we deduce that
	\[ 4 e^{-1}e^{-1} \partial e + e_{(-1)} e^{-1}_{(-1)} (\NO{e h^{-1}})  = (\NO{e h^{-1}}). \]
	By applying $(e_{(-1)})^{-1}$, we get the identity stated in the lemma.
\end{proof}

Now we can finish the proof of \Cref{proposition:relations-localisation-sl2}.

\begin{proof}[Computation of the OPE between $e^{-1}$ and $f$]
	Because $e^{-1}$ is of degree $0$ and $f$ is of degree $2$, only the products $e^{-1}_{(0)}f$ and $e^{-1}_{(1)}f$ may not be zero. We have $e_{(0)} f = h$ and $e_{(1)} f = k 1$, all other products are zero.
	
	Localisation imposes the equality:
	\[ e_{(1)} e^{-1}_{(-1)} f + e_{(0)} e^{-1}_{(0)} f + e_{(-1)} e^{-1}_{(1)} f = 0. \]
	Direct computations give:
	\[ \begin{aligned}
		e_{(1)} e^{-1}_{(-1)} f & = e^{-1}_{(-1)} e_{(1)} f = k e^{-1}, \\
		e_{(0)} e^{-1}_{(0)} f & = e^{-1}_{(0)} e_{(0)} f = e^{-1}_{(0)} h = 2 e^{-1}.
	\end{aligned}  \]
	So, we get the equality
	\[ k e^{-1} + 2 e^{-1} + e_{(-1)} e^{-1}_{(1)} f = 0, \]
	that is to say
	\[ e^{-1}_{(1)} f  = - (k + 2) (e^{-1})^2, \]
	since $(e_{(-1)})^{-1} e^{-1} = (e^{-1})^2$.
	
	Localisation also imposes the identity
	\[ e_{(1)} e^{-1}_{(-2)} f + e_{(0)} e^{-1}_{(-1)} f + e_{(-1)} e^{-1}_{(0)} f + e_{(-2)} e^{-1}_{(1)} f = 0. \]
	By similar computations as before, we get:
	\[ e^{-1}_{(0)} f = (2 k + 2)(e^{-1})^3 \partial e - (e_{(-1)})^{-1} (\NO{e^{-1} h}). \]
	But \Cref{lemma:inverse-no} gives
	\[ \begin{aligned}
		e^{-1}_{(0)} f & = (2 k + 2) (e^{-1})^3 \partial e - (\NO{e^{-1} e^{-1} h} + 4 (e^{-1})^3 \partial e) \\
		& = 2 (k - 1) (e^{-1})^3 \partial e - \NO{e^{-1} e^{-1} h}.
	\end{aligned} \]
\end{proof}

The relation provided by \Cref{proposition:relations-localisation-sl2} can be implemented in Mathematica by using the package \emph{OPEdefs} \cite{thielemans1991mathematica}.

\subsection{The affine W-algebra is a vertex subalgebra of the localisation}

We recall the following well-known fact, see for example \cite[Section 15.3.2]{frenkel2004vertex}.

\begin{theorem} \label{theorem:w-algebra-sl2}
	If $k = -2$, the W-algebra $W^{-2}(\sl_2, f)$ is isomorphic to the centre $\operatorname{Z}(\V^k(\sl_2))$ of the universal affine vertex algebra. This centre is freely generated by
	\[ \widehat{c} \coloneqq \frac{1}{2} \NO{h h} + 2 \NO{e f} - \partial h. \]
	
	If $k \neq -2$, then the W-algebra is isomorphic to the Virasoro vertex algebra of central charge 
	\[ c(k) \coloneqq 1 - \frac{6 (k + 1)^2}{k + 2}. \]
	That means that $\W^k(\sl_2, f)$ is freely generated by one element $\omega$ satisfying the OPE
	\begin{equation} \label{equation:virasoro} \omega(z_1) \omega(z_2) \sim \frac{1}{z_1 - z_2} \partial \omega(z_2) + \frac{2}{(z_1 - z_2)^2} \omega(z_2) + \frac{1}{(z_1 - z_2)^4} \frac{c(k)}{2} \id.  \end{equation}
\end{theorem}

\begin{proposition} \label{proposition:invariant-generator-sl2}
	The invariant subalgebra $\big(\V^k(\sl_2)[e^{-1}]\big)^{\J_\infty(\C h \oplus \C e)}$ is freely generated by the element
	\[ \begin{multlined}
		\widehat{c}_k \coloneqq \frac{1}{2} \NO{h h} + 2 \NO{e f} + (1 + k) \partial h \\ 
		+ (k + 2) \Big(\frac{3 k - 4}{2} (e^{-1})^2 (\partial e)^2 -  \NO{e^{-1} (\partial e) h} - (1+k) \NO{e^{-1} \partial^2 e}\Big).
	\end{multlined}  \]
	\begin{itemize}
		\item If $k = -2$, then $\widehat{c}_{-2} = \widehat{c}$ is the strong generator of the centre $\operatorname{Z}(\V^k(\sl_2))$.
		
		\item If $k \neq -2$, then $\omega \coloneqq \frac{1}{2(k + 2)} \widehat{c}_k$ satisfies the OPE \eqref{equation:virasoro}.
	\end{itemize}
\end{proposition}

\begin{proof}
	The vertex algebra grading and the eigenvalues of the endomorphism $h_{(0)}$ define a bi-grading:
	\[ \begin{multlined}
		\V^k(\sl_2)[e^{-1}] = \bigoplus_{\substack{\Delta \in \Z_{\geqslant 0} \\ \delta \in \Z}} \V^k(\sl_2)[e^{-1}](\Delta, \delta), \\
		\text{where} \quad \V^k(\sl_2)[e^{-1}](\Delta, \delta) \coloneqq \{ a \in \V^k(\sl_2)[e^{-1}](\Delta) \quad | \quad h_{(0)} a = \delta a \}.
	\end{multlined}  \]
	Each piece is finite-dimensional. 
	
	\Cref{theorem:affine-w-algebra-coset} and \Cref{theorem:w-algebra-sl2} state that $ \big(\V^k(\sl_2)[e^{-1}]\big)^{\J_\infty(\C h \oplus \C e)}$ is freely generated by an element $\widehat{c}_k $ that must be in degree $2$ for the vertex algebra grading. This element is in degree $0$ for the action of $h$ because it is $h_{(0)}$-invariant by definition. As a vector space,
	\[ \V^k(\sl_2)[e^{-1}](2, 0) = \operatorname{Span}\big(\NO{h h},  \NO{e f}, \NO{\partial h}, (e^{-1})^2 (\partial e)^2, \NO{e^{-1} (\partial e) h}, \NO{e^{-1} \partial^2 e}\big)  \]
	is $6$-dimensional.
	
	The subspace 
	\[ \Big(\V^k(\sl_2)[e^{-1}](2, 0)\Big)^{\J_\infty(\C h \oplus \C e)}  \] 
	can be computed by brute force computation. Taking a generic element in $\V^k(\sl_2)[e^{-1}]$ and imposing the commutation with $h$ and $e$ gives a set of linear equations (by using \emph{OPEdefs}). Solving this equation leads to a one dimensional vector space
	\[ \Big(\V^k(\sl_2)[e^{-1}](2, 0)\Big)^{\J_\infty(\C h \oplus \C e)}  = \C \widehat{c}_k, \]
	where $\widehat{c}_k$ is given by the formula in \Cref{proposition:invariant-generator-sl2}.
	
	If $k = -2$, the statement is clear. For $k \neq -2$, we can write
	\[ \omega_1 = \omega + \omega_2, \]
	where
	\[ \begin{aligned}
		\omega_1 & \coloneqq \frac{1}{2(k + 2)} \Big(\frac{1}{2} \NO{h h} + 2 \NO{e f} - \partial h\Big) + \frac{1}{2} \partial h, \\
		\omega_2 & \coloneqq - \Big(\frac{3 k - 4}{2} (e^{-1})^2 (\partial e)^2 -  \NO{e^{-1} (\partial e) h} - (1+k) \NO{e^{-1} \partial^2 e}\Big).
	\end{aligned} \]
	By construction of $\omega$ as an element commuting with $e$ and $h$, we get that
	\[ \omega(z_1) \omega_2(z_2) \sim 0. \]
	
	According to \cite[Example 1.4]{kac2003quantum}, $\omega_1$ is the $h$-twisted Sugawara conformal vector in $\V^k(\sl_2)$, meaning that we have the OPE:
	\[ \omega_1(z_1) \omega_1(z_2) \sim \frac{1}{z_1 - z_2} \partial \omega_1(z_2) + \frac{2}{(z_1 - z_2)^2} \omega_1(z_2) + \frac{1}{(z_1 - z_2)^4} \frac{1}{2} \frac{(-3 k)(2 k + 3)}{k + 2} \id. \] 
	A brute force computation with \emph{OPEdefs} gives 
	\[ \omega_2(z_1) \omega_2(z_2) \sim \frac{1}{z_1 - z_2} \partial \omega_2(z_2) + \frac{2}{(z_1 - z_2)^2} \omega_2(z_2) + \frac{1}{(z_1 - z_2)^4} \frac{1}{2} 2 \id. \]
	
	For $n \geqslant 0$, we have $\omega_{(n)} \omega_2 = \omega_2{}_{(n)} \omega = 0$. Therefore
	\[ \omega_1{}_{(n)} \omega_1 = \omega_{(n)} \omega + \omega_2{}_{(n)} \omega_2. \]
	Therefore we get the OPE
	\[ \omega(z_1) \omega(z_2) \sim \frac{1}{z_1 - z_2} \partial \omega(z_2) + \frac{2}{(z_1 - z_2)^2} \omega_1(z_2) + \frac{1}{(z_1 - z_2)^4} \frac{1}{2} c(k) \id, \]
	where
	\[ c(k) = \frac{(-3 k)(2 k + 3)}{k + 2} - 2 = 1 - \frac{6 (k + 1)^2}{k + 2}. \]
\end{proof}

\Cref{theorem:decomposition-chiral} gives an isomorphism
\[ \widehat{\delta} \colon  \V^k(\C h \oplus \C e) \otimes  \big(\V^k(\sl_2)[e^{-1}]\big)^{\J_\infty(\C h \oplus \C e)} \overset{\sim}{\longrightarrow} \big(\V^k(\sl_2)[e^{-1}]\big)^{\J_\infty(\C h \oplus \C e)}. \]
In this case, we can retrieve the inverse reduction \eqref{equation:inverse-reduction} by computing $\widehat{\delta}^{-1}(f)$.

\begin{proposition} \label{proposition:computation-f}
	The following equality
	\[ \begin{multlined}
		f = \NO{e^{-1} \frac{1}{2} \widehat{c}_k} + \frac{(k - 2)(3 k -4)}{4} (e^{-1})^3 (\partial e)^2 \\ + \frac{k - 2}{2} \NO{e^{-1} e^{-1} (\partial e) h} + \frac{(k^2 + k - 4)}{2} \NO{e^{-1} h h} - \frac{k + 1}{2} \NO{e^{-1} \partial h}.
	\end{multlined} \]
	gives an expression of $\widehat{\delta}^{-1}(f)$ in term of $e$, $e^{-1}$, $h$ and $\widehat{c}_k$. 
\end{proposition}

\begin{proof}
	The Borcherds identity \eqref{equation:borcherds} gives:
	\[ f = \NO{e^{-1} e f} + (3k - 4) (e^{-1})^3 (\partial e)^2 - 2 \NO{e^{-1} e^{-1} (\partial e) h} - (k + 3) (e^{-1})^2 \partial^2 e. \]
	Therefore:
	\[ \begin{multlined}
		f = \NO{e^{-1} \frac{1}{2} \widehat{c}_k}  - \NO{e^{-1} \frac{1}{2} \bigg(\frac{1}{2} \NO{h h} + (1 + k) \partial h + \\ (k + 2) \Big(\frac{3 k - 4}{2} (e^{-1})^2 (\partial e)^2  -  \NO{e^{-1} (\partial e) h} - (1+k) \NO{e^{-1} \partial^2 e}\Big)\bigg)}. \\  
		+ (3k - 4) (e^{-1})^3 (\partial e)^2 - 2 \NO{e^{-1} e^{-1} (\partial e) h} - (k + 3) (e^{-1})^2 \partial^2 e. 
	\end{multlined} \]
	After simplification, we get the desired formula.
	
\end{proof}

\printbibliography

\end{document}